%% file: Complete_solutions_of_Toda_equations_and_cyclic_Higgs_bundles_over_non-compact_surfaces.tex
\documentclass{article}
\usepackage{amsmath,amssymb,amscd}

\input{notation}

\input{new_theorem}

\newtheorem{claim}[thm]{Claim}

\title{Complete solutions of Toda equations and cyclic Higgs bundles over non-compact surfaces}
\author{Qiongling Li\thanks{Chern Institute of Mathematics and LPMC, Nankai University, Tianjin 300071, China, qiongling.li@nankai.edu.cn}
\and Takuro Mochizuki\thanks{Research Institute for Mathematical Sciences, Kyoto University, Kyoto 606-8512, Japan, takuro@kurims.kyoto-u.ac.jp}}
\date{}
\begin{document}
\maketitle

\begin{abstract}
On a Riemann surface with a holomorphic $r$-differential, one can naturally define a Toda equation and a cyclic Higgs bundle with a grading. A solution of the Toda equation is equivalent to a harmonic metric of the Higgs bundle for which the grading is orthogonal. Here we focus on a general non-compact Riemann surface with an $r$-differential which is not necessarily meromorphic at infinity. We introduce the notion of complete solution of the Toda equation, and we prove the existence and uniqueness of a complete solution by using techniques for both Toda equations and harmonic bundles. Moreover, we show some quantitative estimates of the complete solution. 
\end{abstract}

\tableofcontents

\section{Introduction}
\subsection{Higgs bundles associated to $r$-differentials
   and harmonic metrics}
 
Let $X$ be any Riemann surface.
We fix a line bundle $K_X^{1/2}$
with an isomorphism
$K_X^{1/2}\otimes K_X^{1/2}\simeq K_X$.
Let $r$ be a positive integer.
We set $\hyperk_{X,r}:=\bigoplus_{i=1}^r K_X^{(r+1-2i)/2}$.
We define the actions of $G_r=\{a\in\cnum\,|\,a^r=1\}$
on $K_{X}^{(r+1-2i)/2}$
by $a\bullet v=a^{i}v$.
They induce a $G_r$-action on $\hyperk_{X,r}$.
For any $r$-differential $q\in H^0(X, K_X^r)$,
let $\theta(q)$ denote the Higgs field
of $\hyperk_{X,r}$ induced by $\theta_i=id: K_X^{(r+1-2i)/2}\rightarrow K_X^{(r+1-2(i+1))/2}\otimes K_X (i=1,\cdots, r-1)$ and $\theta_r=q: K_X^{(1-r)/2}\rightarrow K_X^{(r-1)/2}\otimes K_X$. 

A Hermitian metric $h$ on a Higgs bundle is called harmonic if it satisfies the Hitchin self-dual equation. Geometrically, a harmonic metric gives rise to equivariant harmonic maps from the universal cover $\widetilde X$ to $SL(r,\mathbb C)/SU(r)$. Let $\Harm(q)$ denote the set of $G_r$-invariant
harmonic metrics $h$ of
$(\hyperk_{X,r},\theta(q))$
such that $\det(h)=1$. 
By the $G_r$-invariance,
the decomposition
$\hyperk_{X,r}=\bigoplus_{i=1}^rK_X^{(r+1-2i)/2}$
is orthogonal with respect to
any $h\in \Harm(q)$,
and hence
we obtain the decomposition
$h=\bigoplus h_{|K_X^{(r+1-2i)/2}}$.  
We say that $h\in\Harm(q)$ is \textit{real}
if
$h_{|K_X^{(r+1-2i)/2}}$
and
$h_{|K_X^{(-r-1+2i)/2}}$
are mutually dual.
Let $\Harm^{\real}(q)$ denote the subset of $h\in\Harm(q)$ which are real.

Recall that for a compact Riemann surface $X$, the classification of harmonic metric of $(\hyperk_{X,r},\theta(q))$ is well-known. If $X$ is hyperbolic, there uniquely exists a harmonic metric of unit determinant by the Hitchin-Kobayashi correspondence for Higgs bundles (\cite{Hitchin87, s1}). Moreover, as observed by Baraglia \cite{Baraglia}, the harmonic metric is $G_r$-invariant and real. In other words, if $X$ is compact hyperbolic, $\Harm(q)=\Harm^{\real}(q)$ consists of a unique element. If $X$ is an elliptic curve, it is easy to see that $\Harm(q)$ consists of a unique element if $q\neq 0$, and that $\Harm(0)$ is empty. If $X$ is $\mathbb P^1$, it is easy to see that there exists no non-zero holomorphic $r$-differential, and that $\Harm(0)$ is empty. 

When $X$ is non-compact, the uniqueness of harmonic metrics no longer holds always and the harmonic metric of unit determinant is not necessarily $G_r$-invariant. When $X=\mathbb C^*$ and $q=z^mdz^r$, the solution space $\Toda(q,g)$ was obtained in \cite{GuestItsLin, GuestLin, Toda-lattice, Toda-latticeII} motivated by the relation with the $tt^*$-geometry \cite{CV}. 

In this paper, we study the existence and uniqueness of $G_r$-invariant harmonic metrics on $(\hyperk_{X,r},\theta(q))$ over a non-compact Riemann surface $X$.  

\subsection{Toda equation associated to $r$-differentials}
 
Let $g$ be any K\"{a}hler metric of $X$.
It induces a $G_r$-invariant Hermitian metric $h^{(0)}(g)$
of $\hyperk_{X,r}$.
For any other $G_r$-invariant Hermitian metric $h$
such that $\det(h)=1$,
we obtain a tuple of $\real$-valued functions
$\vecw=(w_1,\ldots,w_r)$ such that $\sum w_i=0$
by the relation
$h_{|K_X^{(r+1-2i)/2}}
=e^{w_i}h^{(0)}(g)_{|K_X^{(r+1-2i)/2}}$.
Then, $h$ is contained in $\Harm(q)$ 
if and only if
the Toda equation is satisfied:
\begin{equation}
\label{eq;20.7.2.1}
\left\{
\begin{array}{l}
 \triangle_gw_1=e^{-w_r+w_1}|q|_g^2-e^{-w_1+w_2}
 -\frac{r-1}{4}k_g
  \\
 \triangle_gw_i=e^{-w_{i-1}+w_i}-e^{-w_i+w_{i+1}}
 -\frac{r+1-2i}{4}k_g \quad (i=2,\ldots,r-1)\\
\triangle_gw_r=e^{-w_{r-1}+w_r}-e^{-w_r+w_1}|q|_g^2
 -\frac{1-r}{4}k_g
\end{array}
\right.
\end{equation}
Here,
$\triangle_g=\frac{1}{2}\sqrt{-1}\Lambda\del\delbar$, 
$k_g=\sqrt{-1}\Lambda R(g)$ and $|q|_g^2=q\bar q/g^r$.
Let $\Toda(q,g)$ denote the set of solutions
$\vecw$ of
(\ref{eq;20.7.2.1})
satisfying $\sum w_i=0$.
A solution $\vecw\in\Toda(q,g)$ is called \textit{real}
if $w_i+w_{r+1-i}=0$. Let $\Toda^{\real}(q,g)$ denote the set of
real solutions of (\ref{eq;20.7.2.1}).
As explained, there is a natural bijection
$\Harm(q)\simeq\Toda(q,g)$, which induces $\Harm^{\real}(q)\simeq\Toda^{\real}(q,g)$.

\begin{rem}
If we change the sign in the system (\ref{eq;20.7.2.1}), we obtain the classical Toda equation studied extensively in integrable system, e.g. see \cite{BPW, BP}. Geometrically, the classical Toda equation gives rise to harmonic maps from surface to compact flag manifolds. 
\end{rem}

\begin{rem}\label{GeometricMeaning}
For general $r$, a solution of the Toda equation gives rise to an equivariant harmonic map $f: \widetilde  X\rightarrow SL(r,\mathbb C)/SU(r)$ such that $tr(\partial f^{\otimes i})=0$ except for $tr(\partial f^{\otimes r})=q$. In lower rank, the Toda equation is encoded with much richer geometry. 
If $r=2$, the Toda equation coincides with the Bochner equation for harmonic maps between surfaces, which is studied extensively in literature, e.g., see \cite{Gupta, HTTW, SchoenYau, WAN, AuWAN, Wolf, Wolf91b}. If $r=3$, the Toda equation for a real solution coincides with Wang's equation for hyperbolic affine spheres in $\mathbb R^3$, which is studied extensively in literature, e.g., see \cite{BH, BH1,DumasWolf, Labourie, Loftin0,  Loftin1, Loftin2, Nie,Wang}.  If $r=4$, the Toda equation for a real solution coincides with the Gauss-Ricci equation for maximal surfaces in $\mathbb H^{2,2}$ by the work of Collier-Tholozan-Toulisse \cite{CTT}, and also studied in \cite{TW}.
\end{rem}

\subsection{Existence and uniqueness of complete solutions}
A solution $\vecw\in\Toda(q,g)$ is called \textit{complete} if for each $2\leq i\leq r$, $e^{-w_{i-1}+w_i}\cdot g$ is complete. In terms of harmonic metrics,
it is equivalent to the condition that
the metrics
$h_{K_X^{(r+1-2(i+1))/2}}\otimes
 h_{K_X^{(r+1-2i)/2}}^{-1}$
induce complete distances on $X$.
Note that
$K_X^{(r+1-2(i+1))/2}\otimes
 (K_X^{(r+1-2i)/2})^{-1}$
 is naturally identified with
 the tangent bundle of $X$. 
 Note that for any K\"ahler metric $g_a$ ($a=1, 2$) there exists a natural bijection $\Toda(q, g_1)\simeq \Toda(q, g_2)$ under which complete solutions are preserved. 

As in Remark \ref{GeometricMeaning}, a complete solution has rich geometric interpretations. For general $r$, the induced metric of the harmonic map $f:\widetilde X\rightarrow SL(r,\mathbb C)/SU(r)$ arising from a complete solution is complete. But our condition of complete solution is stronger than the condition of the induced metric being complete. When $r=2$, a complete solution in $\Toda(q,g)$ is equivalent to looking for an equivariant harmonic map $f: \widetilde X\rightarrow \mathbb H^2$ with Hopf differential $q$ such that the holomorphic energy density $|\partial f|^2$ defines a complete metric on $X$, see Wan \cite{WAN}. This is our main motivation to introduce a complete solution of the Toda equation for general $r$. When $r=3$, a complete real solution in $\Toda(q,g)$ is equivalent to looking for a complete hyperbolic affine sphere in $\mathbb R^3$ with Pick differential $q$. When $r=4$, a complete real solution in $\Toda(q,g)$ gives rise to a complete maximal surface in $\mathbb H^{2,2}$ which is studied recently in \cite{LTW, LT}.

The main result in this paper is the following existence and uniqueness theorem of a complete solution. In the case $r=2$, the theorem is proven in \cite{Li, WAN, AuWAN} in which case the Toda equation reduces to a scalar equation. The existence result makes use of the techniques in rank $2$ case together with the method of super-subsolution for system developed in Guest and Lin \cite{GuestLin}. The uniqueness result makes use of the Omori-Yau and Cheng-Yau maximum principles together with Simpson's inequality for harmonic bundles. 
\begin{thm}[Theorem \ref{ExistenceUniquenessTheorem}]\label{ExistenceCompleteReal}
Suppose $q\neq 0$ unless X is hyperbolic. Then, there exists a unique
complete solution in $\Toda(g, q)$. Moreover, it is real.

If X is parabolic or elliptic, and if $q=0$, then $\Toda(q, g)$ is empty.
 \hfill\qed
\end{thm}

Let $\vecw^{c}$ denote the complete solution and $h^{c}$ the corresponding
harmonic metric.

\subsection{Uniqueness and non-uniqueness of general solutions}
Suppose the $r$-differential $q$ has finitely many zeros, and let $K$ be a relatively compact open subset containing all zeros of $q$. Set $|q|^{2/r}=(q\bar q)^{1/r}$. Then $|q|^{2/r}$ naturally induces a metric on $X\setminus K$. We study the uniqueness of solutions in $\Toda(q,g)$ depending on whether $|q|^{2/r}$ induces a complete metric on $X\setminus K$. In the case $r=2$ or $3$, this is proven by Li \cite{Li} in the setting of a scalar equation for $X=\mathbb C$ and $q$ is a polynomial $r$-differential. 
\begin{thm}[Corollary \ref{UniqueForComplereq}]\label{CompleteDifferentialUnique}
If $|q|^{\frac{2}{r}}$ induces a complete metric on $X\setminus K$, then $\Toda(q,g)=\{\vecw^{c}\}$.
 \hfill\qed
\end{thm}

As a direct corollary of Theorem \ref{CompleteDifferentialUnique}, 
\begin{cor}\label{MeromorphicUnique}
Suppose that $X$ is the complement of a finite subset in a compact
Riemann surface $\overline X$. We also assume that (i) $q$ is meromorphic on $\overline X$,
(ii) the pole order of $q$ is larger than $r$ at each point of $\overline X \setminus X$. Then,
there uniquely exists a solution of the associated Toda equation. Moreover, it is complete and real.
\end{cor}

The statement of Corollary \ref{MeromorphicUnique} also follows from the
method of Kobayashi-Hitchin correspondence for wild harmonic bundles,
as explained in the subsequent paper \cite{LiMochizuki}.

\begin{thm}[Proposition \ref{FiniteManyZerosExistence}]\label{Theorem2:TwoSolutions}
If $|q|^{\frac{2}{r}}$ induces an incomplete metric on $X\setminus K$, then there exists a solution $\vecu\in \Toda^{\real}(q,g)$ which is not complete and satisfies there exists a constant $c>0$ such that on $X\setminus K$, 
\[-\frac{r+1-2l}{r}\log |q|_g-c\leq u_l\leq  -\frac{r+1-2l}{r}\log |q|_g,\quad 1\leq l\leq n=[r/2].\]
In particular, the solutions in $\Toda(q,g)$ are not unique. 
 \hfill\qed
\end{thm}

\subsection{Properties of complete solutions}
We obtain some interesting properties of the complete solution, especially, the precise bounds between the metrics $e^{w_{i+1}-w_i}g(i=1,\cdots, r-1)$. Such estimates were first proven in Dai and Li \cite{DaiLi, DaiLi2} over compact hyperbolic surfaces. When we study real solutions of (\ref{eq;20.7.2.1}), we set $n:=[r/2].$

\begin{thm}[Theorem \ref{CompletenessRealCurvature}]\label{Theorem1:CompleteReal}
The unique complete solution $\vecw^{c}=(w_1,w_2,\cdots, w_r)$ satisfies that $e^{w_{i+1}-w_i}g(i=1,\cdots, r-1)$ are mutually bounded and $|q|_{e^{w_{i+1}-w_i}g}(i=1,\cdots, r-1)$ are bounded. More precisely, one of the following holds:\\
(i) \begin{eqnarray*}
&&w_k<-\frac{r+1-2k}{r}\log|q|_g,\quad 1\leq k\leq  n\\
&&e^{2w_1}|q|_g^2/e^{-w_1+w_2}< 1, \\
&& \frac{(k-1)(r-k+1)}{k(r-k)}<e^{-w_{k-1}+w_k}/e^{-w_k+w_{k+1}}< 1, \quad 2\leq k\leq  n
\end{eqnarray*} 
(ii) $w_k=-\frac{r+1-2k}{r}\log|q|_g$ for $1\leq k\leq n$, in which case $q$ has no zeros and $|q|^{\frac{2}{r}}$ defines a complete metric;\\
(iii) $w_k=\log(\frac{(k-1)!}{(r-k)!}2^{(r+1-2k)})$ for $1\leq k\leq n$, in which case $q\equiv 0$ and $(X, g)$ is a complete hyperbolic surface.

 \hfill\qed
\end{thm}
Let $f:\widetilde X\rightarrow SL(r,\mathbb C)/SU(r)$
be the equivariant harmonic map associated to
the harmonic bundle $(\hyperk_{X,r}, \theta(q), h^{c})$.
Theorem \ref{Theorem1:CompleteReal} implies
that the curvature $\kappa$ of the pullback metric of $f$
satisfies $\kappa<0$ (Corollary \ref{NegativeCurvature})
in Case (i) and (iii).

Next we see that the complete solution dominates any other real solution. 
\begin{thm}[Proposition \ref{RealCompletenessUniqueness}]
Let $\vecw^{c}$ be the complete solution in $\Toda(q,g)$. Suppose $\vecu$ is a real solution in $\Toda(q,g)$. Then 
 one of the following holds: (i) $\vecw^{c}=\vecu$,
 (ii) $w^{\rc}_i<u_i$ $(i=1, \cdots, n)$.
 \hfill\qed
\end{thm}

When the $r$-differential $q$ is bounded with respect to
a complete hyperbolic metric $g$,
we obtain a unique bounded solution.

\begin{thm}[Corollary \ref{BoundedBounded}]\label{BoundedIntro} 
Let $(X, g)$ be a complete hyperbolic surface.
If an $r$-differential $q$ on $X$ is bounded with respect to $g$,
then there uniquely exists a bounded solution in $\Toda(q,g)$.
 Moreover, it is real.
 Conversely, if there exists a bounded solution in $\Toda(q,g)$,
 then $q$ is bounded with respect to $g$. 

In fact, the bounded solution is the complete solution. 
 \hfill\qed
\end{thm}
This result is closely related to a recent work of
Labourie and Toulisse \cite{LT}
on the maximal surfaces in $\mathbb H^{2,2}$.
We also generalize the existence and uniqueness
of bounded solutions to general cyclic Higgs bundles
in a subsequent paper \cite[Section 7]{LiMochizuki}.

In a subsequent paper \cite{LiMochizuki},
we will study a classification of general solutions
of (\ref{eq;20.7.2.1}) for isolated singularities
which are poles or some type of essential singularities.

\subsection{Organization of the paper} In Section \ref{TodaEquation}, we review the relation between Toda equation and the harmonic bundle. In Section \ref{UniquenessSection}, we prove the uniqueness of complete solutions for the Toda equation. For the uniqueness, we will apply the Omori-Yau maximum principle and Cheng-Yau maximum principle. In Section \ref{EstimateSection}, we show some interesting properties of the complete solution. In Section \ref{ExistenceSection}, we show the existence of a complete real solution for the Toda equation. For the existence, we use the method of super-subsolution. 
\subsection{Acknowledgement} 
This study starts from a discussion during the workshop ``Higgs bundles and related topics'' in University of Nice, in 2017.
The authors thank Fran\c{c}ois Labourie and J\'er\'emy Toulisse
for asking a question related to Theorem \ref{BoundedIntro}
and sending their recent preprint.
The first author acknowledges the support from Nankai Zhide foundation. 
The second author is grateful to Martin Guest and Claus Hertling
for discussions on harmonic bundles and Toda equations.
The second author is partially supported by
the Grant-in-Aid for Scientific Research (S) (No. 17H06127),
the Grant-in-Aid for Scientific Research (S) (No. 16H06335),
the Grant-in-Aid for Scientific Research (C) (No. 15K04843),
and the Grant-in-Aid for Scientific Research (C) (No. 20K03609),
Japan Society for the Promotion of Science.


\section{Toda equations and harmonic bundles}\label{TodaEquation}

Let $X$ be a Riemann surface with an $r$-differential $q=q(z)dz^r$.
Let $g=g_0dz\otimes d\bar z$ be a K\"ahler metric of $X$. The K\"ahler form $\omega$ is given by $\omega=\frac{\sqrt{-1}}{2}g_0dz\wedge d\bar z$. 
We also view the K\"ahler metric $g$ as a Hermitian metric on $K_X^{-1}$. 
Let $\Lambda$ denote the contraction with respect to the K\"ahler form $\omega$.
Let $R(g)$ denote the curvature of the Chern connection of $K_X^{-1}$ with $g$. So $R(g)=\bar\partial\partial \log(g_0)=-\partial\bar\partial\log g_0$.
The Laplacian with respect to $g$, the Gaussian curvature of $\Re(g)$, and the square norm of $q$ with respect to $g$ are respectively:
\[\triangle_g:=\frac{\sqrt{-1}}{2}\Lambda\del\delbar=\frac{1}{g_0}\partial_z\partial_{\bar z},\quad k_g:=\sqrt{-1}\Lambda R(g)=-\frac{2}{g_0}\partial_z\partial_{\bar z}\log g_0,\quad |q|_g^2=q(z)\bar q(z)/g_0^r.\]

We obtain the Hermitian metrics $(g^{-1})^{\otimes (r+1-2i)/2}$ on $K_X^{(r+1-2i)/2}$, and $h^{(0)}(g)=\oplus_{i=1}^r(g^{-1})^{\otimes (r+1-2i)/2}$ on $\hyperk_{X,r}$. Then \[|(dz)^{(r+1-2i)/2)}|_{h^{(0)}(g)}=g_0^{-(r+1-2i)/2}.\] For any $\real^r$-valued function
$\vecw=(w_1,\ldots,w_r)$,
we obtain a Hermitian metric
$$h(g,\vecw):=\bigoplus e^{w_i}(g^{-1})^{\otimes (r+1-2i)/2}$$
on $\hyperk_{X,r}$.
\begin{prop}
$h(g,\vecw)\in\Harm(q)$
if and only if
$\vecw$ satisfies the following system of differential equations:
 \begin{equation}
\label{eq;20.9.21.1}
\left\{
\begin{array}{l}
\triangle_g w_1
 =|q|_g^2e^{w_1-w_r}-e^{w_2-w_1}
 -\frac{r-1}{4}k_g \\
  \triangle_g w_j
=e^{w_j-w_{j-1}}-e^{w_{j+1}-w_{j}}
-\frac{r+1-2j}{4}k_g
\quad\quad
(j=2,\ldots,r-1)
  \\
\triangle_g w_r=
 e^{w_r-w_{r-1}}
 -|q|_g^2e^{w_1-w_r}
 +\frac{r-1}{4}k_g 
 \end{array}
 \right.
 \end{equation}
 and $\sum_{i=1}^r w_i=0.$
 
 In particular, in the case $r=2$, $h(g,\vecw)\in\Harm(q)$ if and only if $\vecw=(w_1,-w_1)$ satisfies 
 $$\triangle_g w_1
 =|q|_g^2e^{2w_1}-e^{-2w_1}
 -\frac{1}{4}k_g.$$
\end{prop}
\pf
Let $h(g,\vecw)_i$ denote the restriction of
$h(g,\vecw)$ to $K^{(r+1-2i)/2}$.
We have the decomposition
$h(g,\vecw)=\bigoplus_{i=1}^r h(g,\vecw)_i$.
We obtain the decomposition
of the curvature of the Chern connection:
\[
 R(h(g,\vecw))
 =\bigoplus_{i=1}^r
 R(h(g,\vecw)_i)
 =\bigoplus_{i=1}^r
 \bigl(
 R(g)^{-(r+1-2i)/2})
 +\delbar\del w_i
 \bigr)
\]
Hence, we obtain
\[
 \sqrt{-1}\Lambda R(h(g,\vecw))
 =\bigoplus_{i=1}^r
 \Bigl(
 -\frac{(r+1-2i)}{2}
 k_g
 -2\triangle_gw_i
 \Bigr).
\]
For $P\in X$,
let $v$ denote the element of $(K_X)_{|P}$
such that $|v|_g=1$ by viewing $g^{-1}$ a Hermitian metric on $K_X$. We may choose $v=g_0^{\frac{1}{2}}dz$. Then $\sqrt{-1}\Lambda(v\otimes \bar v)=2$.

We set $v_i:=v^{-(r+1-2i)/2}$.
We obtain $\beta\in\cnum$
such that $q_{|P}=\beta v^r$.
We have
\[
 \theta(v_i)=
 \left\{
  \begin{array}{ll}
   v_{i+1}\otimes v & (i=1,\ldots,r-1)\\
   \beta v_1\otimes v & (i=r)
  \end{array}
 \right.
\]
We obtain
\[
 \theta^{\dagger}_{h(g,\vecw)}(v_i)
 =\left\{
\begin{array}{ll}
 e^{w_i-w_{i-1}}v_{i-1}\otimes \vbar & (i=2,\ldots,r)\\
 \betabar
 e^{w_1-w_r}v_r\otimes\vbar & (i=1)
\end{array}
 \right.
\]
Therefore, we obtain
\[
\sqrt{-1}\Lambda
 [\theta,
 \theta^{\dagger}_{h(g,\vecw)}]
 (v_i)
 =\left\{
\begin{array}{ll}
2 (|\beta|^2e^{w_1-w_r}-e^{w_2-w_1})v_1& (i=1) \\
2 (e^{w_i-w_{i-1}}-e^{w_{i+1}-w_i})v_i & (i=2,\ldots,r-1)\\
2 (e^{w_r-w_{r-1}}-|\beta|^2e^{w_1-w_r})v_{r} & (i=r).
\end{array}
  \right.
\]
Hence, the Hitchin equation
$R(h(g,\vecw))+[\theta,\theta^{\dagger}_{h(g,\vecw)}]=0$
if and only if (\ref{eq;20.9.21.1}) holds.
\hfill\qed

\section{Uniqueness of complete solutions}\label{UniquenessSection}
In this section, we will discuss the uniqueness of solutions for the Toda system (\ref{eq;20.7.2.1}). Our main tools are:
\begin{itemize}
\item the Omori-Yau maximum principle and the Cheng-Yau maximum principle, since they work well for complete Riemann surfaces.
\item  for a harmonic bundle $(E,\delbar_E,\theta,h)$ over a Riemann surface $X$ with a K\"ahler metric $g$, the following inequality holds away from zeros of $\theta$: 
\begin{equation*}
\triangle_g\log|\theta|^2_{h,g}
\geq
\frac{\bigl|[\theta,\theta^{\dagger}_h]\bigr|_{h,g}^2}
{|\theta|_{h,g}^2}+\frac{1}{2}k_g.
\end{equation*}
\end{itemize}
We will apply the two maximum principles to the above inequality. 

\subsection{Maximum principles}
\subsubsection{Omori-Yau Maximum principle}
\begin{lem}[Omori-Yau Maximum Principle \cite{Omori, Yau75}]\label{Omori-Yau}
Suppose $(M,h)$ is a complete manifold with Ricci curvature bounded from below. Then for a $C^2$ function $u:M\rightarrow \mathbb{R}$ bounded from above, there exists a $m_0\in \mathbb{N}$ and a family of points $\{x_m\}\in M$ such that for each $m\geq m_0$,
\[u(x_m)\geq \sup u-\frac{1}{m}, \quad |\nabla_hu(x_m)|\leq \frac{1}{m}, \quad \triangle_h u(x_m)\leq \frac{1}{m},\]
where $\nabla_h$, $\triangle_h$ are the gradient and the Laplacian with respect to the background metric $h$ respectively.
\hfill\qed
\end{lem}

The following is a variant of Omori-Yau maximum principle for manifolds with boundary. 
\begin{lem}
 \label{lem;20.4.10.12} 
 Let $(M,h)$ be any complete Riemannian manifold
with a smooth compact boundary $\del M$.
Assume that the Ricci curvature of $h$ is bounded from below.
Let $u:M\lrarr \real$ be a $C^2$-function bounded from above.
Either one of the following holds.
\begin{itemize}
 \item  $\max_{\del M}u= \sup_{M}u$.
 \item
 There exist  $m_0\in \seisuu_{>0}$ and
 a family of points $Q_m\in M$ $(m\geq m_0)$
 such that
 \[
 u(Q_m)\geq \sup_{M}u -\frac{1}{m},
 \quad
 \bigl|
 \nabla_hu(Q_m)
 \bigr|\leq \frac{1}{m},
 \quad
 \triangle_h u(Q_m)\leq \frac{1}{m},
 \]
 where $\nabla_h$ and $\triangle_h$
 are the gradient and the Laplacian
 with respect to the background metric $h$,
 respectively.
\end{itemize} 
\end{lem}
\pf
Suppose that $\max_{\del M}u<\sup_Mu$.
There exists a compact neighbourhood $K$ of $\del M$ in $M$
with a diffeomorphism
$\Psi:K\simeq \del M\times[0,1]$
such that $\Psi(\del M)=\del M\times\{0\}$.
We may assume that
$\max_K u<\sup_Mu$.
There exists a complete metric $h_1$
on $M\setminus \del M$
such that
$h_{1|M\setminus K}=h_{|M\setminus K}$
and that the Ricci curvature of $h_1$ is
bounded from below.
By applying Omori-Yau maximum principle (Lemma \ref{Omori-Yau}),
we obtain the claim of the lemma.
\hfill\qed

\subsubsection{Cheng-Yau maximum principle}

\begin{lem}[Cheng-Yau Maximum Principle \cite{ChengYau}]\label{Cheng-Yau}
Suppose $(M,h)$ is a complete manifold with Ricci curvature bounded from below. Let $u$ be a $C^2$-function defined on $M$ such that 
$\triangle_h u\geq f(u),$ where $f:\mathbb{R}\rightarrow \mathbb{R}$ is a function. Suppose there is a continuous positive function $g(t):[a,\infty)\rightarrow \mathbb{R}_+$ such that \\
(i) $g$ is non-decreasing; \\(ii) $\lim\inf_{t\rightarrow \infty}\frac{f(t)}{g(t)}>0$;\\ 
(iii) $\int_a^{\infty}(\int_b^tg(\tau)d\tau)^{-\frac{1}{2}}dt<\infty$, for some $b\geq a$,\\
 then the function $u$ is bounded from above. Moreover, if $f$ is lower semi-continuous, $f(\sup u)\leq 0$.
\hfill\qed
\end{lem}
In particular, for $\alpha>1$, positive constants $c_0,c_1,c_2$, one can check $f(t)=c_0t^{\alpha}-(c_1t+c_2)$, $g(t)=t^{\alpha}$ satisfy the above three conditions (i)(ii)(iii). \\

The following is a variant of Cheng-Yau maximum principle for manifolds with boundary. 
\begin{lem}
\label{lem;20.9.24.1} 
Let $(M,h)$ be a connected complete Riemannian manifold
with compact smooth boundary $\del M$
whose Ricci curvature is bounded from below.
Let $u$ be a $C^2$-function $M\lrarr\real$
such that
$\triangle_hu\geq f(u)$
for a function $f$.
Suppose that the restriction of $f$ to $\{-\infty<t<a\}$
is bounded from above for any $a<\infty$.
Suppose that there exists a continuous positive function
$g:\{a\leq t<\infty\}\lrarr\real_{>0}$
such that
\begin{itemize}
 \item $g$ is non-decreasing;
 \item $\liminf_{t\to\infty}\frac{f(t)}{g(t)}>0$;
 \item $\int_a^{\infty}
       \Bigl(\int_b^tg(\tau)d\tau\Bigr)^{-1/2}\,dt
       <\infty$
       for some $b>a$.
\end{itemize}
Then, $u$ is bounded from above.
Moreover, if $f$ is lower semi-continuous,
 either one of the following holds;
 (i) $\max_{\del M}u=\sup_Mu$,
 (ii) $f(\sup_Mu)\leq 0$. 
\end{lem}
\pf
Let $N_1$ be a relatively compact open neighbourhood of $\del M$ in $M$.
Let $N_2$ be a relatively compact open
neighbourhood of $\del M$ in $N_1$.
There exists a complete Riemannian metric $h_0$
of $M\setminus\del M$
such that
(i) Ricci curvature of $h_0$ is bounded from below,
(ii) $h=h_0$ on $M\setminus N_1$.
There exists 
a $C^2$-function $u_0:M\lrarr\real$
such that
(i) $u=u_0$ on $M\setminus N_1$,
(ii) $u_{0|N_2}$ is constant.
The condition (ii) particularly implies that
$\triangle_{h_0}(u_{0|N_2\setminus\del M})=0$.
Because
$\triangle_{h_0}(u_{0})|_{N_1\setminus\del M}$ is bounded,
and because $f(u_0)_{|N_1}$ is bounded from above,
there exists $C_0>0$
such that
$\triangle_{h_0}u_0\geq f(u_0)-C_0$ on $N_1$.
Because $u=u_0$ and $h=h_0$ on $M\setminus N_1$,
we obtain
$\triangle_{h_0}u_0\geq f(u_0)-C_0$ on $M\setminus\del M$.
By Cheng-Yau maximum principle,
$u_0$ is bounded from above.
Because $u_{|N_1}$ is bounded from above,
and because $u=u_0$ on $M\setminus N_1$,
we obtain that $u$ is bounded from above.

Set $T:=\sup_M u$.
Assume that $T>\max_{\del M}u$.
We may assume that
$T>\sup_{N_1}u$.
There exists $\epsilon>0$
such that
$\sup_{N_1}u<T-2\epsilon$.
Let $\chi:\real\lrarr\real$  be a $C^{\infty}$-function
such that
(i) $0\leq \chi(t)\leq 1$ for any $t$,
(ii) $\chi(t)=0$ $(t\geq T-\epsilon)$,
(iii) $\chi(t)=1$ $(t\leq T-2\epsilon)$.
We set $f_0(t):=f(t)-C_0\chi(t)$.
On $M\setminus N_1$,
we obtain
\[
 \triangle_{h_0}(u_0)
 =\triangle_h(u)
 \geq f(u)\geq f(u_0)-\chi(u_0)C_0
 =f_0(u_0).
\]
On $N_1$,
we obtain
\[
 \triangle_{h_0}(u_0)
 \geq f(u_0)-C_0
 =f(u_0)-\chi(u_0)C_0
 =f_0(u_0).
\]
By Cheng-Yau maximum principle (Lemma \ref{Cheng-Yau}),
we obtain $f_0(\sup_{M\setminus\del M} u_0)\leq 0$.
Because $\sup_{M\setminus\del M} u_0=T=\sup_M u$,
we obtain $f(\sup_M u)\leq 0$.
\hfill\qed


\subsection{Preliminary from linear algebra}

\subsubsection{Lower estimate of the commutator with the adjoint
(general case)}

We recall well known lemmas from \cite[Page 729]{s2}
for the convenience of the readers.
Let $V$ be an $r$-dimensional vector space
with a Hermitian metric $h$.
Let $f$ be an endomorphism of $V$.
Let $f^{\dagger}_h$ denote the adjoint of $f$
with respect to $h$.

There exists a (non-canonical) orthonormal base
$v_1,\ldots,v_r$ of $V$
for which $f$ is lower triangular,
i.e.,
$f(v_j)=\sum_{j\leq i} f_{i,j}v_i$.
We define the endomorphism
$f_0$ of $V$
by $f_0(v_i)=f_{i,i}v_i$ $(i=1,\ldots,r)$.
We set $f_1:=f-f_0$.
Note that $|f|_h^2=|f_0|^2+|f_1|_h^2$.

\begin{lem}
\label{lem;20.9.18.6}
There exists $C_0>0$ depending only on $r$
such that
$\bigl|
 \bigl[f,f^{\dagger}_h\bigr]
 \big|_h\geq C_0|f_1|_{h}^2$.
\end{lem}
\pf
We study the claim in terms of matrices.
The following argument is explained in \cite[Page 729]{s2}.
Let $A=(A_{i,j})$ be an $r$-square matrix
such that $A_{i,j}=0$ $(i<j)$.
Let $A^{\dagger}$ denote the adjoint matrix of $A$,
i.e.,
$A^{\dagger}_{i,j}=\overline{(A_{j,i})}$.
We decompose $A=A_0+A_1$,
\[
 (A_0)_{i,j}=\left\{
\begin{array}{ll}
 A_{i,i}& (i=j) \\
 0 & (i\neq j),
\end{array}
\right.
\quad\quad
 (A_1)_{i,j}=\left\{
\begin{array}{ll}
 0& (i=j) \\
 A_{i,j} & (i\neq j).
\end{array}
\right.
\]
We obtain
$[A,A^{\dagger}]
=[A_0,A_1^{\dagger}]
+[A_1,A_0^{\dagger}]
+[A_1,A_1^{\dagger}]$.
Note that
the diagonal entries of
$[A_0,A_1^{\dagger}]$
and $[A_1,A_0^{\dagger}]$ are $0$.
By direct calculations,
we obtain the following
for any $i=1,\ldots,r$:
\[
 [A_1,A_1^{\dagger}]_{i,i}
=\sum_{1\leq j<i}|A_{i,j}|^2
-\sum_{i<j\leq r} |A_{j,i}|^2.
\]
Note that the following holds
for any $i$:
\begin{multline}
\label{eq;20.9.18.4}
 \sum_{i<j\leq r}
 |A_{j,i}|^2
\leq
 \Bigl|
 \sum_{1\leq j<i}
 |A_{i,j}|^2
-\sum_{i<j\leq r}
 |A_{j,i}|^2
 \Bigr|
+ \sum_{1\leq j<i}
 |A_{i,j}|^2
 \\
\leq \Bigl|
 \sum_{1\leq j<i}
 |A_{i,j}|^2
-\sum_{i<j\leq r}
 |A_{j,i}|^2
 \Bigr|
 +\sum_{1\leq j<i}
 \sum_{j<k\leq r}
 |A_{k,j}|^2
\end{multline}
By an easy induction using (\ref{eq;20.9.18.4}),
we can prove that
there exist positive constants $C_{1,k}$ $(k=1,\ldots,r)$
depending only on $r$
such that the following holds:
\begin{equation}
\label{eq;20.9.18.5}
 \sum_{i\leq k}
 \sum_{i<j\leq r}
 |A_{j,i}|^2
\leq
C_{1,k}
 \sum_{i\leq k}
 \left|
 \sum_{1\leq j< i}|A_{i,j}|^2
-\sum_{i<j\leq r}|A_{j,i}|^2
 \right|.
\end{equation}
The inequality in the case $k=r$
implies the claim of Lemma \ref{lem;20.9.18.6}.
\hfill\qed

\begin{lem}
\label{lem;20.9.18.20}
 Suppose that there exists $B>0$
such that 
any eigenvalues $\alpha$ of $f$
satisfies $|\alpha|<B$,
and that $|f|_h^2\geq 2rB^2$.
Then, we obtain 
\[
 \bigl|
 [f,f^{\dagger}_h]
 \bigr|_h
 \geq \frac{C_0}{2}|f|_h^2.
\]
Here, $C_0$ is a constant as in
Lemma {\rm\ref{lem;20.9.18.6}}.
\end{lem}
\pf
Because $|f_0|_h^2\leq rB^2$,
we obtain
$|f_1|_h^2\geq rB^2\geq |f_0|_h^2$.
It implies that
$|f_1|_h^2\geq \frac{1}{2}|f|_h^2$.
Thus, we obtain the claim of the lemma.
\hfill\qed

\subsubsection{Lower estimate of the commutator with
 the adjoint in the cyclic case}

Let $e_1,\ldots,e_r$ be an orthogonal base of $(V,h)$,
i.e., $h(e_i,e_j)=0$ $(i\neq j)$.
Let $\alpha$ be a nonzero complex number.
Let $f$ be the endomorphism of $V$
determined by
$f(e_i)=e_{i+1}$ $(i=1,\ldots,r-1)$
and $f(e_r)=\alpha ^re_1$.

\begin{prop}
\label{prop;20.9.17.40}
For any $0<\epsilon<1$,
there exists $\delta>0$,
depending only on $\epsilon$ and $r$,
such that the following holds.
\begin{itemize}
 \item If there exists $1\leq i\leq r-1$ such that
       $|e_{i+1}|_h\cdot |e_i|_h^{-1}\leq \epsilon|\alpha|$,
       then we obtain
       $\bigl|
       [f,f^{\dagger}_h]
       \bigr|\geq \delta|f|_h^2$.
\end{itemize}
\end{prop}
\pf
We begin with a preliminary.
We set
$Z:=
\bigl\{\veca=(a_1,\ldots,a_r)\in\real^r_{>0}\,\big|\,
\prod_{i=1}^r a_i=1
\bigr\}$.
We define the function
$F:Z\lrarr \real_{>0}$
by
\[
 F(\veca)=\sum_{i=1}^{r-1} (a_{i+1}-a_i)^2
 +(a_1-a_r)^2.
\]
We also define the function $G:Z\lrarr \real_{>0}$
by
\[
 G(\veca)=\sum_{i=1}^ra_i.
\]
\begin{lem}
\label{lem;20.9.22.10}
 For any $0<\epsilon<1$,
there exists $\delta>0$ such that the following holds
for $\veca\in Z$:
\begin{itemize}
 \item If there exists $1\leq i\leq r$
       such that $a_i<\epsilon$,
       then $F(\veca)>\delta G(\veca)^2$.
\end{itemize}
\end{lem}
\pf
To simplify the notation,
we set $a_0:=a_r$ and $a_{r+1}:=a_1$.
We set $c_{\max}:=\max\{a_i\}$
and $c_{\min}:=\min\{a_i\}$.
We obtain $G(\veca)^2\leq r^2c_{\max}^2$.
There exists $i$
such that $|a_{i+1}-a_i|\geq r^{-1}(c_{\max}-c_{\min})$.
Hence, we obtain
\[
F(\veca)\geq r^{-2}(c_{\max}-c_{\min})^2
=r^{-2}c_{\max}^2(1-c_{\max}^{-1}c_{\min})^2.
\]
Because $\prod_{i=1}^ra_i=1$,
we obtain
$1\leq c_{\min}c_{\max}^{r-1}$,
i.e.,
\[
c_{\max}^{-1}\leq
c_{\min}^{\frac{1}{r-1}}\leq \epsilon^{\frac{1}{r-1}}.
\]
We obtain
$c_{\max}^{-1}c_{\min}
\leq
\epsilon^{\frac{r}{r-1}}$.
Then, the claim of the lemma is clear.
\hfill\qed

\vspace{.1in}
Let us return to the proof of Proposition \ref{prop;20.9.17.40}.
Note that
\[
 |f|_h^2
=\sum_{i=1}^{r-1}
 \frac{|e_{i+1}|_h^2}{|e_i|_h^2}
+\frac{|\alpha|^{2r}|e_1|_h^2}{|e_r|^2}.
\]
By a direct calculation, we obtain
\[
 f^{\dagger}_h(e_{i+1})=\frac{|e_{i+1}|_h^2}{|e_i|_h^2}e_i
 \quad(i=1,\ldots,r-1),
 \quad
 f^{\dagger}_h(e_1)
 =\frac{\alphabar^r|e_1|_h^2}{|e_r|_h^2}e_r.
\]
Hence, we obtain
\[
 (f\circ f^{\dagger}_h-f^{\dagger}_h\circ f)(e_1)
 =\left(
 \frac{|e_1|_h^2|\alpha|^{2r}}{|e_r|_h^2}
-\frac{|e_2|_h^2}{|e_1|_h^2}
 \right)e_1,
\]
\[
 (f\circ f^{\dagger}_h-f^{\dagger}_h\circ f)(e_i)
 =\left(
 \frac{|e_i|_h^2}{|e_{i-1}|_h^2}
-\frac{|e_{i+1}|_h^2}{|e_i|_h^2}
 \right)e_i
 \quad(i=2,\ldots,r-1),
\]
\[
 (f\circ f^{\dagger}_h-f^{\dagger}_h\circ f)(e_r)
 =\left(
 \frac{|e_r|_h^2}{|e_{r-1}|_h^2}
-\frac{|e_{1}|_h^2|\alpha|^{2r}}{|e_r|_h^2}
 \right)e_r
\]
It implies
\[
 \bigl|
 [f,f^{\dagger}_h]
 \bigr|_h^2
 =
\left(
 \frac{|e_1|_h^2|\alpha|^{2r}}{|e_r|_h^2}
-\frac{|e_2|_h^2}{|e_1|_h^2}
 \right)^2
+\sum_{i=2}^{r-1}
 \left(
 \frac{|e_i|_h^2}{|e_{i-1}|_h^2}
-\frac{|e_{i+1}|_h^2}{|e_i|_h^2}
 \right)^2
+\left(
 \frac{|e_r|_h^2}{|e_{r-1}|_h^2}
-\frac{|e_{1}|_h^2|\alpha|^{2r}}{|e_r|_h^2}
 \right)^2
\]
Then, we obtain the claim of Proposition \ref{prop;20.9.17.40}
from Lemma \ref{lem;20.9.22.10}.
\hfill\qed

\subsubsection{Commutator with the difference of metrics}
Let $h'$ be another Hermitian metric of $V$
such that $e_1,\ldots,e_r$ is orthogonal with respect to $h'$.
Assume that
$\prod |e_i|_h=\prod|e_i|_{h'}=1$.
We obtain the automorphism $s$
determined by $h'(x,y)=h(sx,y)$
for $x,y\in V$,
which satisfies $\det(s)=1$.
Let $B>0$ and $C>0$.
Assume the following conditions on $h$ and $h'$:
\begin{itemize}
 \item $|e_i|_{h}\cdot |e_{i+1}|_{h}^{-1}\leq B$
       for $i=1,\ldots,r-1$
       \item $|s|_{h}\leq C$.
\end{itemize}

\begin{lem}
 \label{lem;20.7.5.30}
 There exist
 $C_1>0$ and $\epsilon_1>0$,
 depending only on $B$ and $C$,
 such that
 the following holds:
 \begin{itemize}
  \item  If there exists $0<\epsilon<\epsilon_1$
	 such that
	 $\bigl|[f,s]s^{-1/2}\bigr|_{h}\leq \epsilon$,
	 then $\bigl|s-\id_V\bigr|_{h}\leq
	 C_1\cdot\epsilon$ holds.
 \end{itemize}
\end{lem}
\pf
There exist $a_i\in\real_{>0}$ $(i=1,\ldots,r)$
such that
$s(e_i)=a_ie_i$
and that $\prod a_i=1$.
Note that $a_i<C$ and $a_i^{-1}<C^{r-1}$.
For $i=1,\ldots,r-1$,
we obtain
\[
 [f,s]s^{-1/2}e_i
 =a_i^{1/2}(1-a_i^{-1}a_{i+1})e_{i+1}.
\]
Hence, we obtain the following for any $i=1,\ldots,r-1$:
\[
 |1-a_i^{-1}a_{i+1}|
 \leq
  (|e_{i}|_{h}\cdot |e_{i+1}|_{h}^{-1})
 a_i^{-1/2}
 \bigl|
 [f,s(h)]s(h)^{-1/2}
 \bigr|_{h}
\leq BC^{(r-1)/2}\epsilon.
\]
There exists $C_2>0$
such that
if $BC^{(r-1)/2}\epsilon<1/2$
then the following holds:
\[
\bigl|
-\log(a_i)+\log(a_{i+1})
\bigr|
=
\bigl|
\log(a_i^{-1}a_{i+1})
\bigr|
\leq
 C_2|1-a_i^{-1}a_{i+1}|
 \leq
  BC^{(r-1)/2}C_2\epsilon.
\]
We obtain
$|\log(a_i)-\log(a_1)|
\leq
 iBC^{(r-1)/2}C_2\epsilon$.
Because $\prod a_i=1$,
we obtain
\[
 r|\log a_1|
 \leq
 \frac{r(r-1)}{2}BC^{(r-1)/2}C_2\epsilon.
\]
We obtain
$|\log(a_i)|\leq
 \frac{r+1}{2}BC^{(r-1)/2}C_2\epsilon$
 for $i=1,\ldots,r$.
It implies the claim of the lemma.
\hfill\qed

\subsection{Boundedness of Higgs field
  with bounded spectral curve}

Let $(E,\delbar_E,\theta,h)$ be a harmonic bundle on
a Riemann surface $X$.
Let $g$ be a K\"ahler metric on $X$.
Let $Z(\theta)$ denote the zero set of $\theta$.
Recall that there is the following standard inequality
on $X\setminus Z(\theta)$ 
(see \cite[Page 729]{s2}):
\begin{equation}
\label{eq;20.9.23.1}
\triangle_g\log|\theta|^2_{h,g}
\geq
\frac{\bigl|[\theta,\theta^{\dagger}_h]\bigr|_{h,g}^2}
{|\theta|_{h,g}^2}+\frac{1}{2}k_g.
\end{equation}
Here,
$\triangle_g=\frac{1}{2}\sqrt{-1}\Lambda_g\del\delbar$,
and
$k_g=\sqrt{-1}\Lambda_gR(g)$.

\begin{condition}
\label{condition;20.9.16.1}
There exists $C_0>0$ such that
any eigenvalue $\alpha$ of $\theta$ at $P\in X$
satisfies
\[
|\alpha|_{g}\leq C_0.
\]
More precisely,
for any holomorphic coordinate $z$ around $P$,
we express $\theta$ as $\theta=f\,dz$ around $P$,
and then any eigenvalue $\alpha_1$  of $f$ at $P$
satisfies
$|\alpha_1|\cdot |dz|_g\leq C_0$.
\hfill\qed  
\end{condition}

\begin{lem} 
\label{lem;20.9.17.10}
 If Condition {\rm\ref{condition;20.9.16.1}} is satisfied,
there exist $C_i>0$ $(i=1,2)$,
depending only on $C_0$ and $\rank(E)$,
such that 
the following inequality holds on $X\setminus Z(\theta)$:
\[
 \triangle_g\log|\theta|_{h,g}^2
 \geq
 C_1|\theta|_{h,g}^2
-C_2
 +\frac{1}{2}k_g.
\]
 \end{lem}
\pf
By Lemma \ref{lem;20.9.18.20},
Condition {\rm\ref{condition;20.9.16.1}} implies that
there exists $C_3>0$
depending only on $\rank(E)$,
such that the following holds.
\begin{itemize}
 \item If $|\theta|_{h,g}^2\geq 2(\rank E)C^2_0$ at $P\in X$,
       then
       $\bigl|[\theta,\theta^{\dagger}_h]\bigr|_{h,g}
      \geq C_3|\theta|_{h,g}^2$ at $P$.
\end{itemize}
Therefore,
we obtain the following on $X\setminus Z(\theta)$:
\[
 \bigl|[\theta,\theta^{\dagger}_h]\bigr|_{h,g}^2
 \geq
 C^2_3|\theta|_{h,g}^2
 \bigl(|\theta|_{h,g}^2-2C_0^2(\rank E)\bigr).
\]
Then, we obtain the claim of the lemma
from (\ref{eq;20.9.23.1}).
\hfill\qed

\vspace{.1in}
Let $X_1\subset X$ be a relatively compact open subset
whose boundary is smooth and compact.
Let $\Xbar_1$ denote the closure of $X_1$ in $X$.

\begin{prop}
\label{prop;20.9.17.11}
In addition to Condition {\rm\ref{condition;20.9.16.1}},
we assume that $g_{|\Xbar_1}$ is complete and that
the Gaussian curvature of $g_{|\Xbar_1}$ is bounded from below.
Then, $|\theta|_{h,g}$ is bounded on $\Xbar_1$.
More precisely,
there exists $C_{10}>0$,
depending only on $\inf_{\Xbar_1}k_g$, $C_0$ and $\rank(E)$
such that
$|\theta|_{h,g}\leq
 \max\bigl\{
  C_{10},\max_{\del X_1}|\theta|_{h,g}
 \bigr\}$
on $X_1$.
\end{prop}
\pf
By the assumptions and Lemma \ref{lem;20.9.17.10},
there exists $C_{11}>0$,
depending only on
$\inf_{\Xbar_1} k_{g}$, $C_0$ and $\rank(E)$
such that
the following holds on $X_1\setminus Z(\theta)$:
\[
 \triangle_g\log|\theta|_{h,g}^2
 \geq
 C_1|\theta|_{h,g}^2
 -C_{11}.
\]
It implies the following inequality on $X_1\setminus Z(\theta)$:
\begin{equation}
\label{eq;20.9.23.10}
 \triangle_g|\theta|_{h,g}^2
  \geq
  |\theta|_{h,g}^2\bigl(
  C_1|\theta|_{h,g}^2
  -C_{11}
  \bigr).
\end{equation}
Because the both sides are continuous on $X_1$,
the inequality (\ref{eq;20.9.23.10}) holds on $X_1$.
By Lemma \ref{lem;20.9.24.1},
$|\theta|_{h,g}$ is bounded on $X_1$.
(We set
$u=|\theta|_{h,g}^2$,
$f(t)=C_1t^2-C_{11}t$
and $g(t)=C_1t^2$ in Lemma \ref{lem;20.9.24.1}.)
If $\sup_{X_1}|\theta|_{h,g}>\max_{\del X_1}|\theta|_{h,g}$,
then 
we obtain
$\sup|\theta|_{h,g}^2\leq C_1^{-1}C_{11}$,
where the right hand side depends only on
$\inf_{\Xbar_1}k_g$, $C_0$ and $\rank(E)$.
\hfill\qed\\

We shall use the following lemma to prove that there is no
solution for the Toda equation on a parabolic Riemann surface with the
trivial $r$-differential.
\begin{lem}\label{NoSolutionOnParabolicSurface}
Suppose that $X$ is parabolic,
i.e., $X$ is $\cnum$, $\cnum^{\ast}$ or an elliptic curve.
If $\theta$ is nilpotent, then $\theta=0.$
\end{lem}
\pf
There exists a complete K\"ahler metric $g$ of $X$
such that $k_g=0$.
Because $\theta$ is nilpotent,
Lemma \ref{lem;20.9.18.20} implies that
there exists $C>0$ such that
the following holds on $Z(\theta)$:
\[
 \triangle_g\log|\theta|_{h,g}^2\geq C|\theta|_{h,g}^2.
\]
As in the proof of Proposition \ref{prop;20.9.17.11},
we obtain
\begin{equation}
\label{eq;20.9.30.1}
 \triangle_g|\theta|_{h,g}^2\geq C|\theta|_{h,g}^4
\end{equation}
on $X\setminus Z(\theta)$.
Because the both side of (\ref{eq;20.9.30.1}) are continuous
on $X$,
the inequality (\ref{eq;20.9.30.1}) holds on $X$.
It implies that 
$|\theta|_{h,g}^2$ is subharmonic on $X$.
By Proposition \ref{prop;20.9.17.11},
$|\theta|_{h,g}^2$ is bounded.
Because $X$ is $\cnum$, $\cnum^{\ast}$, or an elliptic curve,
we obtain that $|\theta|_{h,g}^2$ is constant.
By (\ref{eq;20.9.30.1}),
we obtain $\theta=0$.
\hfill\qed

\subsection{Toda equation}

Let $X$ be a Riemann surface.
Let $q$ be an $r$-differential on $X$.
Note that $\theta(q)$ is nowhere vanishing.

\subsubsection{The case where the $r$-differential is bounded}

Let $g$ be a complete K\"ahler metric of $X$
such that $|q|_g$ is bounded.

\begin{prop}
\label{prop;20.9.18.1}
Assume that the Gaussian curvature of $g$ is bounded from below.
Then, there exists $C_0>0$
depending only on
$\inf k_g$, $\sup|q|_g$ and $r$
such that the following holds
for any $h\in\Harm(q)$:
\[
 |\theta(q)|_{h,g}\leq C_0.
\] 
As a result,
there exists $C_1>0$
depending only on
$\inf k_g$, $\sup|q|_g$ and $r$
such that the following holds
for any $\vecu\in\Toda(q,g)$:
\[
 u_{i+1}-u_{i}<C_1
\quad (i=1,\ldots,r-1),
\quad
 u_1-u_r+\log|q|^2_g<C_1.
\]
\end{prop}
\pf
By the boundedness of $|q|_g$,
there exists $C_2>0$
such that any eigenvalue $\alpha$ of
$\theta(q)$ at $P\in X$
satisfies
$|\alpha|_g\leq C_2$.
We obtain the first claim by
Proposition \ref{prop;20.9.17.11}.

Recall that
$\theta(q)$ is obtained as
$\theta(q)=\sum_{i=1}^r
 \theta(q)_i$,
 where 
$\theta(q)_i:K_X^{(r+1-2i)/2}\lrarr K_X^{(r+1-2(i+1))/2}\otimes K_X$
$(i=1,\ldots,r-1)$
are induced by the identity,
and
$\theta(q)_r:K_X^{(-r+1)/2}\lrarr K_X^{(r-1)/2}\otimes K_X$
is induced by the multiplication of $q$.
Because
$\hyperk_{X,r}=\bigoplus K_X^{(r+1-2i)/2}$
is orthogonal with respect to $h\in\Harm(q)$,
we obtain
\[
 |\theta(q)|^2_{h,g}
 =\sum_{i=1}^r
 |\theta(q)_i|^2_{h,g}.
\]
Let $\vecu\in\Toda(q,g)$ denote the solution
corresponding to $h$,
i.e.,
$h_{|K_X^{(r+1-2i)/2}}
=e^{u_i}\cdot (g^{-1})^{\otimes(r+1-2i)/2}$.
Note that
\[
 |\theta(q)_i|^2_{h,g}=e^{u_{i+1}-u_i}
 \quad (i=1,\ldots,r-1),
 \quad
 |\theta(q)_r|^2_{h,g}
 =e^{u_1-u_r}|q|^2_{g}
\]
Then, we obtain the second claim of the proposition
from the first claim.
\hfill\qed\\

Let $h\in\Harm(q)$.
For $i=1,\ldots,r-1$,
there exists the natural isomorphism:
\[
\Bigl(
 K_X^{(r+1-2i)/2}
 \Bigr)^{-1}
 \otimes
 K_X^{(r+1-2(i+1))/2}
 \simeq K_X^{-1}.
\]
Let $g(h)_i$ be the K\"ahler metric of $X$
obtained as
$h_{|K_X^{(r+1-2i)/2}}^{-1}
\otimes
h_{|K_X^{(r+1-2(i+1))/2}}$.
In terms of $\vecw\in \Toda(q,g)$, the metrics $g(h)_i=e^{-w_i+w_{i+1}}\cdot g (i=1,\cdots, r-1)$. We will show that  the curvature of the metrics $g(h)_i (i=1,\cdots, r-1)$ are bounded from below. This fact will be useful in applying the two maximum principles.

\begin{lem}\label{LowerBoundCurvature}
For any solution $\vecw\in\Toda(q,g)$, the metrics $e^{-w_1+w_2}\cdot g, \cdots, e^{-w_{r-1}+w_r}\cdot g$ satisfy that their Gaussian curvature are bounded from below by $-4$.
\end{lem}
\pf 
Recall the curvature formula for a K\"ahler metric $g=g_0dz\otimes d\bar z$ is $k_g=-\frac{2}{g_0}\triangle\log g_0$, where $\triangle=\partial_z\partial_{\bar z}$. So for $r\geq 3$,
\begin{eqnarray*}k_{e^{-w_1+w_2}\cdot g}&=&-\frac{2}{e^{-w_1+w_2}\cdot g_0}\triangle[(-w_1+w_2)+\log g_0]=-\frac{2}{e^{-w_1+w_2}}[\triangle_{g}(-w_1+w_2)-\frac{1}{2}k_g]\\
&=&2\frac{e^{-w_r+w_1}|q|_g^2-2e^{-w_1+w_2}+e^{-w_2+w_3}}{e^{-w_1+w_2}}> -4.\\
k_{e^{-w_k+w_{k+1}}\cdot g}&=&-\frac{2}{e^{-w_k+w_{k+1}}\cdot g_0}[\triangle(-w_k+w_{k+1})+\log g_0]\\&=&2\frac{e^{-w_{k-1}+w_k}-2e^{-w_k+w_{k+1}}+e^{-w_{k+1}+w_{k+2}}}
{e^{-w_k+w_{k+1}}}> -4, \quad 2\leq k\leq r-2.\\
k_{e^{-w_{r-1}+w_r}\cdot g}&=&-\frac{2}{e^{-w_{r-1}+w_r}\cdot g_0}\triangle[(-w_{r-1}+w_r)+\log g_0]\\
&=&2\frac{e^{-w_{r-2}+w_{r-1}}-2e^{-w_{r-1}+w_r}+e^{-w_r+w_1}|q|_g^2}{e^{-w_{r-1}+w_r}}> -4.
\end{eqnarray*}
For $r=2$, 
\begin{eqnarray*}k_{e^{-2w}\cdot g}&=&-\frac{2}{e^{-2w}\cdot g_0}\triangle[(-2w)+\log g_0]=-\frac{2}{e^{-2w}}[\triangle_{g}(-2w)-\frac{1}{2}k_g]\\
&=&4\frac{e^{2w_1}|q|_g^2-e^{-2w}}{e^{-2w}}\geq -4.
\end{eqnarray*}
\hfill\qed 

\begin{cor}\label{Comparison}
 \label{cor;20.9.17.12}
Assume that
there exist a solution $\vecw\in\Toda(q,g)$
and $1\leq k\leq r-1$
such that 
$|w_{k+1}-w_{k}|$ is bounded.
Then, there exists $C>0$
such that the following holds
for any $\vecu\in\Toda(q,g)$:
\[
 u_{i+1}-u_{i}<C
\quad (i=1,\ldots,r-1),
 \quad
 u_1-u_r+\log|q|^2_g<C.
\] 
\end{cor}
\pf
The K\"ahler metric
$\gtilde=e^{w_{k+1}-w_k}g$ is complete since $|w_{k+1}-w_k|$ is bounded. By Lemma \ref{LowerBoundCurvature},
the Gaussian curvature of $\gtilde$ is bounded from below by
a constant. And $|q|_{\gtilde}$ is still bounded since $|w_{k+1}-w_k|$ is bounded. 
Hence, we obtain the claim of the corollary
from Proposition \ref{prop;20.9.18.1}.
\hfill\qed

\vspace{.1in}
Let $g^{(j)}$ $(j=1,2)$ be complete K\"ahler metrics of $X$
such that
$|q|_{g^{(j)}}$ are bounded.

\begin{prop}
\label{prop;20.9.19.3}
 Suppose that there exist solutions
 $\vecw^{(j)}\in\Toda(q,g^{(j)})$
 and $1\leq k(j)\leq r-1$
such that
$|w^{(j)}_{k(j)+1}-w^{(j)}_{k(j)}|$
are bounded.
Then, $g^{(1)}$ and $g^{(2)}$ are mutually bounded.
\end{prop}
\pf
We obtain an $\real_{>0}$-valued function
$A$ determined by $g^{(2)}=Ag^{(1)}$.
We set
\[
\wtilde^{(2)}_i:=
w^{(2)}_i-\frac{r+1-2i}{2}\log A.
\]
Then, $\widetilde{\vecw}^{(2)}\in \Toda(q,g^{(1)})$.
By Corollary \ref{cor;20.9.17.12},
there exists $C_1>0$ such that
$\wtilde^{(2)}_{i+1}-\wtilde^{(2)}_{i}<C_1$
for $i=1,\ldots,r-1$.
In particular, we obtain
\[
w^{(2)}_{k(2)+1}-w^{(2)}_{k(2)}
+\log A
=\wtilde^{(2)}_{k(2)+1}-\wtilde^{(2)}_{k(2)}<C_1.
\]
Because $|w^{(2)}_{k(2)+1}-w^{(2)}_{k(2)}|$ is bounded,
there exists $C_2>0$
such that $\log A<C_2$.
By exchanging the roles of
$(g^{(1)},\vecw^{(1)})$
and
$(g^{(2)},\vecw^{(2)})$,
we obtain that
there exists $C_3>0$
such that $-\log A<C_3$.
Therefore, we obtain
that $g^{(1)}$ and $g^{(2)}$ are mutually bounded.
\hfill\qed

\begin{cor}
\label{cor;20.9.18.10}
Suppose that there exist bounded solutions
$\vecw^{(j)}\in\Toda(q,g^{(j)})$.
Then, $g^{(1)}$ and $g^{(2)}$ are mutually bounded.
\hfill\qed
\end{cor}

\subsubsection{The case where 
an associated K\"ahler metric is complete}
\label{subsection;20.9.18.101}

\begin{lem}
\label{lem;20.9.17.50}
There exists $C_j$ $(j=1,2)$,
depending only on $r$,
such that the following holds
on $X$
for any $i=1,\ldots,r-1$:
\[
 \triangle_{g(h)_i}\log|\theta(q)|_{h,g(h)_i}
 \geq
 C_1|\theta(q)|^2_{h,g(h)_i}-C_2.
\]
\end{lem}
\pf
Let $P$ be any point of $X$.
We fix $i$.
Take a holomorphic coordinate $z$ around $P$
such that $|dz|_{g(h)_i}^2=2$ at $P$.
We set $e_k=(dz)^{(r+1-2k)/2}$ $(k=1,\ldots,r)$.
We express
$q$ and $\theta(q)$ as
$q=\beta(dz)^r$
and
$\theta(q)=f\,dz$ around $P$.
We obtain
$f(e_k)=e_{k+1}$ $(k=1,\ldots,r-1)$
and
$f(e_r)=\beta e_1$.
By the construction,
we obtain the following at $P$:
\[
 \frac{|e_{i+1}|^2_h}{|e_i|_h^2}=
 |dz|_{g(h)_i}^{-2}=\frac{1}{2}.
\]
By Proposition \ref{prop;20.9.17.40},
there exist $C_{10}>0$ and $C_{11}>0$,
depending only on $r$,
such that the following holds.
\begin{itemize}
 \item If $|\beta|^{1/r}\geq C_{10}$ at $P$,
       then
       $\bigl|[f,f^{\dagger}_h]\bigr|_h\geq C_{11}|f|_h^2$
       at $P$.
\end{itemize}
By Lemma \ref{lem;20.9.18.20},
there exists $C_{12}>0$,
depending only on $r$,
such that the following holds.
\begin{itemize}
 \item If $|\beta|^{1/r}<C_{10}$
       and $|f|_h^2>2rC^2_{10}$ at $P$,
       then we obtain
       $\bigl|[f,f^{\dagger}_h]\bigr|_h\geq C_{12}|f|_h^2$
       at $P$.
\end{itemize}
By our choice of $C_{12}$
the following holds.
\begin{itemize}
 \item If $|\beta|^{1/r}<C_{10}$ at $P$,
       then we obtain the following at $P$:
       \[
       \bigl|[f,f^{\dagger}_h]\bigr|_h^2
       \geq C_{12}^2|f|_h^2\bigl(|f|_h^2-2rC^2_{10}\bigr).
       \]
\end{itemize}
Recall that $k_{g(h)_i}\geq -4$.
Then, we obtain the claim of Lemma \ref{lem;20.9.17.50}.
\hfill\qed

\begin{prop}
\label{prop;20.9.17.100}
If there exists $1\leq i\leq r-1$
such that $g(h)_i$ is complete,
then $|\theta(q)|_{h,g(h)_i}$ is bounded.
As a result,
the functions $g(h)_{j}/g(h)_i$ $(j=1,\ldots,r-1)$
and $|q|_{g(h)_i}$ are bounded.
\end{prop}
\pf
We obtain the boundedness of
$|\theta(q)|_{h,g(h)_i}$
by Cheng-Yau maximum principle
and Lemma \ref{lem;20.9.17.50}.
Note that
$|\theta(q)_j|^2_{h,g(h)_i}
 =g(h)_j/g(h)_i$.
(See the proof of Proposition \ref{prop;20.9.18.1}
for $\theta(q)_j$.)
Because
\[
 \sum_{j=1}^{r}|\theta(q)_j|^2_{h,g(h)_i}= |\theta(q)|^2_{h,g(h)_i},
\] for each $j=1, \cdots, r$, $|\theta(q)_j|^2_{h,g(h)_i}$ is bounded. For $j=1, \cdots, r-1$, we obtain the boundedness of
$g(h)_j/g(h)_i$.

Let $g$ be any K\"ahler metric of $X$.
Let $\vecu\in\Toda(q,g)$ be the solution corresponding to $h$.
Also we have 
\begin{equation*}
e^{u_1-u_r}|q|_g^2, e^{-u_1+u_2}, \cdots, e^{-u_{r-1}+u_r}\leq e^{-u_i+u_{i+1}}\cdot C. 
\end{equation*}
Thus, $|q|_g^2\leq (e^{-u_i+u_{i+1}})^r\cdot C^r$ and $|q|_{g(h)_i}$ is bounded 
\hfill\qed

\begin{cor}
\label{cor;20.9.17.101}
If $g(h)_i$ $(i=1,\ldots,r-1)$ are complete,
then $g(h)_i$ are mutually bounded,
and $|q|_{g(h)_i}$ are bounded.
In other words,
there exists a complete K\"ahler metric $g$
 such that $g$ and $g(h)_i$ are mutually bounded,
 and that $|q|_g$ is bounded.
\end{cor}
\pf
By Proposition \ref{prop;20.9.17.100},
$g(h)_i$ are mutually bounded.
Then, the claims of the corollary are clear.
\hfill\qed

\begin{cor}\label{BoundedSolutionImplyBoundedDifferential}
Let $g$ be a complete K\"ahler metric of $X$.
If there exists a bounded solution $\vecw\in\Toda(q,g)$,
then $|q|_g$ is bounded.
\hfill\qed
\end{cor}

\subsubsection{Uniqueness of complete solutions}

\begin{thm}
\label{thm;20.9.18.100}
Suppose that the following holds for
 $h^{(a)}\in\Harm(q)$ $(a=1,2)$.
\begin{itemize}
 \item $g(h^{(a)})_i$ $(i=1,\ldots,r-1)$ are complete.
\end{itemize}
Then, we obtain $h^{(1)}=h^{(2)}$.
\end{thm}
\pf 
By Corollary \ref{cor;20.9.17.101},
there exist complete K\"ahler metrics
$g^{(a)}$ $(a=1,2)$ such that
(i) $g^{(a)}$ and $g(h^{(a)})_i$ $(i=1,\ldots,r-1)$
are mutually bounded,
(ii) $|q|_{g^{(a)}}$ are bounded.
By Corollary \ref{cor;20.9.18.10},
we obtain that $g^{(1)}$ and $g^{(2)}$
are mutually bounded.
We may assume $g^{(1)}=g^{(2)}=:g$.
Moreover,
we may assume that the Gaussian curvature of $g$
is bounded from below.
Note that $|q|_g$ is bounded.

Let $\vecw^{(a)}\in\Toda(q,g)$
be the solutions corresponding to $h^{(a)}$.
Because $g(h^{(a)})_i$ $(i=1,\ldots,r-1)$
are mutually bounded with $g$,
we obtain that
$|w^{(a)}_{i+1}-w^{(a)}_i|$ $(i=1,\ldots,r-1)$
are bounded.
Because $\sum w^{(a)}_i=0$,
we obtain that $\vecw^{(a)}$ are bounded.
It implies that
$h^{(1)}$ and $h^{(2)}$ are mutually bounded.

We apply the Omori-Yau maximum principle
to prove $h^{(1)}=h^{(2)}$
from the mutually boundedness of $h^{(1)}$
and $h^{(2)}$ as follows.

Let $s$ be the automorphism of $\hyperk_{X,r}$
determined by
$h^{(2)}=h^{(1)}\cdot s$.
By \cite[Lemma 3.1]{s1}, there exists a constant $C>0$ such that following 
 inequality holds on $X$:
\begin{equation}
 \label{eq;20.7.5.22}
 \triangle_g\Tr(s)
  =C\bigl|
  [\theta(q),s]s^{-1/2}
  \bigr|^2_{h^{(1)},g}
  +C\bigl|
   \delbar_{\hyperk_{X,r}}(s)s^{-1/2}
   \bigr|^2_{h^{(1)},g}
   \geq
   C\bigl|[\theta(q),s]s^{-1/2}
   \bigr|^2_{h^{(1)},g}.
\end{equation}

By Omori-Yau maximum principle,
there exists
a sequence 
$Q_{\ell}\in X$ $(\ell=1,2,\ldots)$
such that
\[\triangle_g\Tr(s)(Q_{\ell})\leq \ell^{-1},
\quad
\Tr(s)(Q_{\ell})>\sup_{Q\in X}\Tr(s)(Q)-\ell^{-1}.
\]
Hence, there exists $C_1>0$ such that
the following holds for any $\ell$:
\[
 \bigl|
 [s,\theta(q)]s^{-1/2}
 \bigr|^2_{h^{(1)},g}(Q_{\ell})
 \leq C_{1}\ell^{-1}.
\]
Let $z_{\ell}$ be a holomorphic coordinate
around $Q_{\ell}$
such that $|dz_{\ell}|^2_{g}(Q_{\ell})=2$.
Because $g(h^{(1)})_i$ and $g$ are mutually bounded,
there exists $B>0$
such that
$\bigl|(dz_{\ell})^{(r+1-2i)/2}\bigr|_{h^{(1)}}
\cdot
\bigl|(dz_{\ell})^{(r+1-2(i+1))/2}\bigr|_{h^{(1)}}^{-1}
\leq B$
for any $\ell$.
By Lemma \ref{lem;20.7.5.30},
there exists $C_{2}>0$ and $\ell_2$
such that
$\Tr(s)(Q_{\ell})\leq r(1+C_2\ell^{-1/2})$
for any $\ell>\ell_2$.
Hence, we obtain
$\sup_{Q\in X}\Tr(s)(Q)\leq r$.
Because $\det(s)=1$,
we obtain $\Tr(s)\geq r$ for any $Q$.
Therefore, $\Tr(s)$ is constantly $r$,
and we obtain $s=\id$.
\hfill\qed

\begin{cor}
\label{cor;20.9.18.102}
 Suppose that $h\in\Harm(q)$ satisfies
 the following condition.
\begin{itemize}
 \item $g(h)_i$ $(i=1,\ldots,r-1)$ are complete.
\end{itemize}
Then, $h$ is real, 
i.e.,
 $h_{|K_X^{(r+1-2i)/2}}$
 and
 $h_{|K_X^{(-(r+2-2i)/2)}}$
 are mutually dual.
\end{cor}
\pf
We obtain $h^{\lor}\in\Harm(q)$
as the dual of $h$
by using the natural identification of
$\hyperk_{X,r}$
with its dual.
Because $h^{\lor}\in\Harm(q)$ also satisfies
the same condition,
we obtain $h=h^{\lor}$
by Theorem \ref{thm;20.9.18.100},
which implies the claim of the corollary.
\hfill\qed

\vspace{.1in}
We can restate
Theorem \ref{thm;20.9.18.100}
and Corollary \ref{cor;20.9.18.102}
in terms of solutions of the Toda equation.

\begin{cor}
\label{cor;20.9.18.110}
Let $g$ be a K\"ahler metric of $X$.
 \begin{itemize}
  \item If $\vecw^{(a)}\in\Toda(q,g)$ $(a=1,2)$ are complete, 
then we obtain
$\vecw^{(1)}=\vecw^{(2)}$. 
  \item
       If there exists a complete solution
$\vecw\in\Toda(q,g)$,
it is real,
i.e.,
       $w_i+w_{r+1-i}=0$.
       \hfill\qed
\end{itemize}
\end{cor}

We also obtain the following result
on bounded solutions.
\begin{cor}\label{BoundedSolutionUnique}
Let $g$ be a complete K\"ahler metric of $X$.
\begin{itemize}
  \item If $\vecw^{(a)}\in\Toda(q,g)$ $(a=1,2)$ are bounded,
then we obtain
$\vecw^{(1)}=\vecw^{(2)}$. 
  \item If there exists a bounded solution
$\vecw\in\Toda(q,g)$,
it is real, i.e., $w_i+w_{r+1-i}=0$.
       \hfill\qed
\end{itemize}
\end{cor}

\subsubsection{Estimates on an open subset
with smooth compact boundary}

Let $X_1\subset X$ be an open subset
whose boundary is smooth and compact.
Let $\Xbar_1$ denote the closure of $X_1$ in $X$.

\begin{prop}
\label{prop;20.9.19.10}
 Suppose that $h\in\Harm(q)$ satisfies the following condition.
\begin{itemize}
 \item $g(h)_{i|\Xbar_1}$  $(i=1,\ldots,r-1)$ are complete.
\end{itemize}
Then, $g(h)_{i|\Xbar_1}$ are mutually bounded,
and $|q|_{g(h)_i}$ are bounded on $\Xbar_1$.
\end{prop}
\pf
We obtain Proposition \ref{prop;20.9.19.10}
as a consequence of the following lemma,
which we can prove
by using the argument in the proof of 
Proposition \ref{prop;20.9.17.100}
with Cheng-Yau maximum principle
(Lemma \ref{lem;20.9.24.1}).
\begin{lem}
\label{lem;20.9.19.11}
Let $h\in\Harm(q)$.
If there exists $1\leq i\leq r-1$
such that
$g(h)_{i|\Xbar_1}$ is complete,
then the functions
$|\theta(q)|_{h,g(h)_i}$,
$|q|_{g(h)_i}$
and $g(h)_j/g(h)_i$ $(j=1,\ldots,r-1)$
are bounded on $\Xbar_1$.
\hfill\qed
\end{lem}

The following proposition
is a classification of complete solutions
up to boundedness.

\begin{prop}
\label{prop;20.9.19.20}
 Let $h^{(a)}\in\Harm(q)$ $(a=1,2)$.
Suppose that
$g(h^{(a)})_{i|\Xbar_1}$ $(i=1,\ldots,r-1)$ are complete.
Then, $h^{(1)}$ and $h^{(2)}$
are mutually bounded on $X_1$.
More precisely,
for the automorphism $s$ of $\hyperk_{X,r}$
determined by $h^{(2)}=h^{(1)}\cdot s$,
we obtain
$\sup_{\Xbar_1}\Tr(s)=\max_{\del X_1}\Tr(s)$.
In particular,
if $h^{(1)}_{|\del X_1}=h^{(2)}_{|\del X_1}$,
we obtain $h^{(1)}=h^{(2)}$ on $X_1$.
\end{prop}
\pf
We obtain the following lemma from
Proposition \ref{prop;20.9.17.11}.

\begin{lem}
\label{lem;20.9.19.2}
 Let $g$ be a K\"ahler metric of $X$ such that
(i) $|q|_{g}$ is bounded on $\Xbar_1$,
(ii) $g_{|\Xbar_1}$ is complete
whose Gaussian curvature is bounded from below.
Then, for any $h\in\Harm(q)$,
$|\theta(q)|_{h,g}$ is bounded on $\Xbar_1$.
Equivalently,
for any $\vecu\in\Toda(q,g)$,
 $u_{i+1}-u_i$ $(i=1,\ldots,r-1)$
 and
 $u_1-u_r+\log|q|_g^2$
 are bounded from above on $\Xbar_1$. 
\hfill\qed
\end{lem}

As a consequence of Lemma \ref{lem;20.9.19.2},
we obtain the following lemma.
\begin{lem}
\label{lem;20.9.19.22}
Let $g$ be a K\"ahler metric of $X$ such that
(i) $|q|_{g}$ is bounded on $\Xbar_1$,
 (ii) $g_{|\Xbar_1}$ is complete.
Suppose that there exist $\vecw\in\Toda(q,g)$ and
$1\leq k\leq r-1$
such that
$|w_{k+1}-w_k|$ is bounded on $\Xbar_1$.
Then, for any $\vecu\in \Toda(q,g)$,
 $u_{i+1}-u_i$ $(i=1,\ldots,r-1)$
 and
 $u_1-u_r+\log|q|^2_g$
 are bounded from above on $\Xbar_1$.
\hfill\qed
\end{lem}

We obtain the following lemma
by using the argument in the proof of
Proposition \ref{prop;20.9.19.3}
with Lemma \ref{lem;20.9.19.22}.
\begin{lem}
\label{lem;20.9.19.21}
 Let $g^{(j)}$ $(j=1,2)$ be K\"ahler metrics of $X$
 such that
 (i) $|q|_{g^{(j)}}$ are bounded on $\Xbar_1$,
 (ii) $g^{(j)}_{\Xbar_1}$ are complete.
Suppose that there exist solutions
$\vecw^{(j)}\in\Toda(q,g^{(j)})$
and $1\leq k(j)\leq r-1$
such that
$|w^{(j)}_{k(j)+1}-w^{(j)}_{k(j)}|_{|\Xbar_1}$
are bounded.
 Then, $g^{(1)}_{|\Xbar_1}$ and $g^{(2)}_{|\Xbar_1}$
 are mutually bounded.
 In particular, if there exist bounded solutions
 $\vecw^{(j)}\in \Toda(q,g^{(j)})$,
 then $g^{(1)}_{|\Xbar_1}$
 and $g^{(2)}_{|\Xbar_1}$
 are mutually bounded.
\hfill\qed
\end{lem}

By Proposition \ref{prop;20.9.19.10}
and Lemma \ref{lem;20.9.19.21},
$h^{(1)}$ and $h^{(2)}$ are mutually bounded.
If $\sup_{\Xbar_1}\Tr(s)>\max_{\del X_1}\Tr(s)$,
Omori-Yau maximum principle (Lemma \ref{lem;20.4.10.12})
implies that there exists a sequence $Q_{\ell}\in X_1$
as in the proof of Theorem \ref{thm;20.9.18.100},
and hence we obtain
$\sup_{\Xbar_1}\Tr(s)=r\leq \max_{\del X_1}\Tr(s)$,
which is a contradiction.
Hence, we obtain 
$\sup_{\Xbar_1}\Tr(s)=\max_{\del X_1}\Tr(s)$.
Thus, the proof of Proposition \ref{prop;20.9.19.20}
is completed.
\hfill\qed

\vspace{.1in}

In terms of the solutions of the Toda equation,
Proposition \ref{prop;20.9.19.10}
and Proposition \ref{prop;20.9.19.20}
are stated as follows.
\begin{cor}
Let $g$ be a K\"ahler  metric of $X$.
\begin{itemize}
 \item  Suppose that $\vecw\in\Toda(q,g)$ satisfies
	the following condition.
	\begin{itemize}
	 \item
	      $(e^{w_{i+1}-w_i}g)_{|\Xbar_1}$ $(i=1,\ldots,r-1)$
	are complete.
	\end{itemize}
	Then, $(e^{w_{i+1}-w_i}g)_{|\Xbar_1}$
	are mutually bounded,
	and
	$e^{w_1-2w_r+w_{r-1}} |q|_g$ is bounded on $\Xbar_1$.
 \item Suppose that $\vecw^{(a)}\in\Toda(q,g)$ $(a=1,2)$
       satisfy the following condition.
       \begin{itemize}
	\item $(e^{w^{(a)}_{i+1}-w^{(a)}_i}g)_{|\Xbar_1}$
	      $(i=1,\ldots,r-1)$
	are complete.
       \end{itemize}
       Then, $|w^{(1)}_i-w^{(2)}_i|$ are bounded on $\Xbar_1$.
       More precisely,
       we obtain
       \[
	\sup_{\Xbar_1}\sum_{i=1}^r e^{w^{(1)}_i-w^{(2)}_i}
       =
       \max_{\del X_1}\sum_{i=1}^r e^{w^{(1)}_i-w^{(2)}_i}.
       \]
       \hfill\qed
\end{itemize}
\end{cor}

\subsubsection{The case where $|q|^{2/r}$ is complete
outside a relatively compact open subset}

Let $X_1\subset X$ be an open subset
whose boundary is smooth and compact.
\begin{prop}
\label{prop;20.9.20.1}
Suppose that $|q|^{2/r}_{|\Xbar_1}$ is complete.
\begin{itemize}
 \item Let $g$ be any K\"ahler metric of $X$
       such that $g=|q|^{2/r}$ on $\Xbar_1$.
       Then, for any $\vecu\in\Toda(q,g)$,
       $|u_i|$ $(i=1,\ldots,r)$ are bounded on $\Xbar_1$.
 \item For any $h\in\Harm(q)$,
 $g(h)_i$ $(i=1,\ldots,r-1)$ are mutually bounded with
 $|q|^{2/r}$ on $\Xbar_1$.
 \item Any $h^{(a)}\in\Harm(q)$ $(a=1,2)$ are mutually bounded
on $\Xbar_1$.
 \end{itemize}
\end{prop}
\pf
Let $g$ be a K\"ahler metric of $X$
such that $g_{|\Xbar_1}=|q|^{2/r}_{|\Xbar_1}$.
Let $\vecu\in\Toda(q,g)$.
By Lemma \ref{lem;20.9.19.2},
there exists $C_1>0$ such that
$u_{i+1}-u_{i}\leq C_1$ $(i=1,\ldots,r-1)$
and
$u_1-u_r+\log|q|^2_g\leq C_1$ on $\Xbar_1$.
Because 
$\bigl|\log|q|^2_g\bigr|$ is bounded on $\Xbar_1$,
there exists $C_2>0$ such that
$u_{i+1}-u_{i}\leq C_2$ $(i=1,\ldots,r-1)$
and
$u_1-u_r\leq C_2$ on $\Xbar_1$.
For any $j=2,\ldots,r$,
we obtain
$u_j\leq u_1+(j-1)C_2$
from $u_{i+1}-u_i\leq C_2$ $(i=1,\ldots,j-1)$
on $\Xbar_1$.
Hence, we obtain
$u_j\leq u_1+rC_2$ $(j=2,\ldots,r)$
on $\Xbar_1$.
By the cyclic symmetry,
for any $i$,
we obtain
$u_j\leq u_i+rC_2$ $(j\neq i)$
on $\Xbar_1$.
Thus, we obtain
$|u_i-u_j|\leq rC_2$ for any $1\leq i,j\leq r$
on $\Xbar_1$.
Because $\sum_{i=1}^r u_i=0$,
there exists $C_3>0$
such that
$|u_i|\leq C_3$ $(i=1,\ldots,r)$
on $\Xbar_1$.
It implies the claims of
Proposition \ref{prop;20.9.20.1}.
\hfill\qed
  
\begin{cor}
Suppose that there exists a relatively compact open subset $N$
such that $|q|^{2/r}$ is complete on $X\setminus N$.
Then, $\Harm(q)$ consists of at most a unique complete solution. 
\end{cor}
\pf
It follows from
Theorem \ref{thm;20.9.18.100} and
Proposition \ref{prop;20.9.20.1}.
\hfill\qed

\begin{cor}\label{UniqueForComplereq}
Suppose that there exists a relatively compact open subset $N$
such that $|q|^{2/r}$ is complete on $X\setminus N$.
For any K\"ahler metric $g$ of $X$,
$\Toda(q,g)$ consists of at most a unique complete solution.
\hfill\qed
\end{cor}

\subsubsection{Pull back}

Let $F:X_1\lrarr X_2$ be a holomorphic map
of Riemann surfaces
which is locally an isomorphism,
i.e., the derivative of $F$ is nowhere vanishing.
Let $q$ be an $r$-differential on $X_2$.

Because $F$ is locally an isomorphism,
there exists a natural isomorphism
$F^{\ast}(\hyperk_{X_2,r},\theta(q))
\simeq
(\hyperk_{X_1,r},\theta(F^{\ast}q))$.
For any $h\in\Harm(q)$,
it is well known and easy to check that
the induced metric $F^{\ast}(h)$
of $\hyperk_{X_1,r}$ is a harmonic metric of 
$(\hyperk_{X_1,r},\theta(F^{\ast}q))$.
Moreover,
$F^{\ast}(h)$ is clearly 
$G_r$-invariant
and satisfies $\det(F^{\ast}(h))=1$,
and hence
$F^{\ast}(h)\in\Harm(F^{\ast}q)$.
In this way,
we obtain the map
\[
F^{\ast}:\Harm(q)\lrarr \Harm(F^{\ast}q).
\]

Let $g$ be a K\"ahler metric of $X_2$.
We obtain a K\"ahler metric
$F^{\ast}(g)$ of $X_1$.
For any $\vecw\in \Toda(q,g)$,
we obtain an $\real^r$-valued function
$F^{\ast}(\vecw)$ on $X_1$.
Because
$F^{\ast}h(g,\vecw)=h(F^{\ast}g,F^{\ast}\vecw)$
as a Hermitian metric on
$\hyperk_{X_1,r}$,
we obtain $F^{\ast}(\vecw)\in \Toda(F^{\ast}q,F^{\ast}g)$.
Thus, we obtain
$F^{\ast}:\Toda(q,g)\lrarr \Toda(F^{\ast}q,F^{\ast}g)$.

If $X_1=X_2$ and $F^{\ast}(q)=q$,
then $F$ induces an automorphism on $\Harm(q)$.
If moreover $F^{\ast}g=g$,
then $F$ induces an automorphism of
$\Toda(q,g)$.

\begin{prop}\label{Pullback}
Suppose that $X_1=X_2$, $F^{\ast}(q)=q$ and $F^{\ast}(g)=g$.
Then, a complete solution $\vecw\in\Toda(q,g)$
is preserved by $F$,
i.e., $F^{\ast}(\vecw)=\vecw$.
\end{prop}
\pf
It follows from the uniqueness of complete solutions.
\hfill\qed

\subsection{Appendix: Tame harmonic bundles
  on parabolic Riemann surfaces}

In this appendix,
after recalling the general theory of tame harmonic bundles
due to Simpson \cite{s2},
we explain the classification of
tame harmonic bundles on parabolic Riemann surfaces
(elliptic curves, $\cnum$, or $\cnum^{\ast}$)
as an easy consequence of Lemma \ref{NoSolutionOnParabolicSurface}.
Though it is easy and more or less well known to specialists,
the clear statements might be convenient somewhere.
  
\subsubsection{Tame harmonic bundles on punctured Riemann surfaces}
\label{subsection;20.9.29.100}

Let $X$ be a compact Riemann surface.
Let $D$ be a finite subset in $X$.
For any $P\in D$,
let $(X_P,z_P)$ be a holomorphic coordinate neighbourhood
around $P$ such that $z_P(P)=0$.
We set $X_P^{\ast}:=X_P\setminus\{P\}$.

Let $(E,\delbar_E,\theta,h)$ be a harmonic bundle on $X\setminus D$.
We obtain the spectral curve
$\Sigma_{\theta}\subset K_{X\setminus D}$.
The harmonic bundle $(E,\delbar_E,\theta,h)$
is called tame \cite{s2}
if the closure of $\Sigma_{\theta}$
in the logarithmic cotangent bundle $K_X(D)$
is proper over $X$.

Let $\nbige^0$ denote the sheaf of holomorphic sections of
$(E,\delbar_E)$.
We recall that $(\nbige^0,\theta)$
extends to a filtered regular Higgs bundle
in \cite{s2}.

For any $a\in\real$,
$\nbige^0_{|X_P^{\ast}}$
extends to an $\nbigo_{X_P}$-module
$\nbigp^h_a\nbige^0$ as follows.
For any open subset $U\subset X_P$ such that $P\in U$,
$\nbigp^h_a\nbige^0(U)$ is the space of
holomorphic sections $s$ of $\nbige^0_{|U\setminus\{P\}}$
satisfying
$|s|_h=O\bigl(|z_P|^{-a-\epsilon}\bigr)$
for any $\epsilon>0$.
According to \cite{s2},
$\nbigp^h_a(\nbige^0_{|X_P^{\ast}})$
is a locally free $\nbigo_{X_P}$-module,
$\theta_{|X_P^{\ast}}$ induces
a logarithmic Higgs field of
$\nbigp^h_a(\nbige^0_{|X_P^{\ast}})$,
i.e.,
$\theta_{|X_P^{\ast}}:
\nbigp^h_a(\nbige^0_{|X_P^{\ast}})
\lrarr
 \nbigp^h_a(\nbige^0_{|X_P^{\ast}})\otimes K_{X_P}(P)$.
There exist natural monomorphisms
$\nbigp^h_b(\nbige_{|X_P^{\ast}})
\lrarr
\nbigp^h_a(\nbige_{|X_P^{\ast}})$
for any $b\leq a$. 
We define
\begin{equation}
\label{eq;20.9.29.10}
 \Gr^{\nbigp^h}_a(\nbige^0)_P:=
 \nbigp^h_a(\nbige^0_{|X_P^{\ast}})
 \big/
 \bigcup_{b<a}\nbigp^h_b(\nbige^0_{|X_P^{\ast}}),
\end{equation}
which we naturally regard as a finite dimensional
$\cnum$-vector space.

For $\veca\in\real^{D}$,
we obtain a locally free $\nbigo_X$-module
$\nbigp^h_{\veca}(\nbige^0)$
from $\nbige^0$
and $\nbigp^h_{a_P}(\nbige^0_{|X_P^{\ast}})$ $(P\in D)$.
Thus, we obtain a regular filtered Higgs bundle
$(\nbigp^h_{\ast}\nbige^0,\theta)$ on $(X,D)$.
As proved in \cite{s2},
we obtain
\[
 \deg(\nbigp^h_{\veca}\nbige^0)
 -\sum_{P\in D}
 \sum_{a_P-1<b\leq a_P}
 b\dim\Gr^{\nbigp^h}_{b}(\nbige^0)_P
 =0.
\]

Recall the polystability of
$(\nbigp^h_{\ast}\nbige^0,\theta)$.
Let $E_1\subset E$ be a holomorphic subbundle
such that $\theta(E_1)\subset E_1\otimes K_{X\setminus D}$.
Let $\nbige^0_1\subset\nbige^0$ be
the sheaf of holomorphic sections of $E_1$.
Suppose the following holds.
\begin{condition}
\label{condition;20.9.29.21}
For any $a\in\real$,
       $\nbige^0_{1|X_P^{\ast}}$ extends to
       a locally free $\nbigo_{X_P}$-submodule
       $\nbigp^h_a(\nbige^0_{1|X_P^{\ast}})
       \subset
       \nbigp^h_a(\nbige^0_{|X_P^{\ast}})$
       such that
       $\nbigp^h_a(\nbige^0_{|X_P^{\ast}})\big/
       \nbigp^h_a(\nbige^0_{1|X_P^{\ast}})$
       is locally free.
       We define
       $\Gr^{\nbigp^h}_a(\nbige^0_1)_P$
       as in {\rm(\ref{eq;20.9.29.10})}.
\end{condition}
For any $\veca\in\real^{D}$,
we obtain a locally free $\nbigo_X$-module
$\nbigp^h_{\veca}(\nbige^0_1)$
from $\nbige^0_1$
and $\nbigp^h_{a_P}(\nbige^0_{1|X_P^{\ast}})$ $(P\in D)$.
We define
\[
\deg(\nbigp^h_{\ast}\nbige^0_1):=
 \deg(\nbigp^h_{\veca}\nbige^0_1)
 -\sum_{P\in D}
 \sum_{a_P-1<b\leq a_P}
 b\dim\Gr^{\nbigp^h}_b(\nbige^0_{1})_P.
\]
Note that the left hand side is well defined
in the sense that it is independent of $\veca$.

\begin{thm}[\cite{s2}]
\label{thm;20.9.30.10}
 For any Higgs subbundle $E_1\subset E$
satisfying Condition {\rm\ref{condition;20.9.29.21}},
we obtain
$\deg(\nbigp^h_{\ast}\nbige^0_1)\leq 0$.
If $\deg(\nbigp^h_{\ast}\nbige^0_1)=0$,
then the orthogonal complement
$E_2$ of $E_1$ in $E$
is a holomorphic subbundle
such that $\theta(E_2)\subset E_2$.
As a result,
we obtain a decomposition of harmonic bundle
$(E,\delbar_E,\theta,h)
=(E_1,\delbar_{E_1},\theta_{1},h_{1})
 \oplus
 (E_2,\delbar_{E_2},\theta_2,h_2)$,
 which induces
 the decomposition of filtered Higgs bundles
 $(\nbigp^h_{\ast}\nbige^0,\theta)
 =(\nbigp^h_{\ast}\nbige^0_1,\theta_1)
 \oplus
 (\nbigp^h_{\ast}\nbige^0_2,\theta_2)$.
 \hfill\qed
\end{thm}


\subsubsection{Harmonic bundles on elliptic curves}

Let $X$ be an elliptic curve.
A harmonic bundle $(E,\delbar_E,\theta,h)$
of rank $1$ on $X$
is easily described.
Indeed,
$(E,\delbar_E)$ is a holomorphic line bundle,
$h$ is a Hermitian metric of $E$
whose Chern connection is flat,
and $\theta$ is a holomorphic one form on $X$.
The following proposition is well known.

\begin{prop}
For any harmonic bundle
$(E,\delbar_E,\theta,h)$ of rank $r$ on $X$,
there exist
 harmonic bundles of rank one
 $(E_i,\delbar_{E_i},\theta_i,h_i)$
 $(i=1,\ldots,r)$
 such that
 $(E,\delbar_E,\theta,h)
 \simeq
 \bigoplus_{i=1}^r
 (E_i,\delbar_{E_i},\theta_i,h_i)$. 
\end{prop}
\pf
There exists a nowhere vanishing holomorphic one form
$dz$ on $X$.
Let $f$ be the endomorphism of $E$
determined by $\theta=f\,dz$.
We obtain the characteristic polynomial
$\det(t\id_E-f)=\sum_{j=0}^r a_jt^j$.
Then, $a_j$ are holomorphic functions on $X$,
and hence they are constant,
which implies that the eigenvalues of $f$
are constant.
Therefore,
we obtain a decomposition of Higgs bundles
$(E,\delbar_E,\theta)=
\bigoplus (E_{\alpha},\delbar_{E_{\alpha}}\theta_{\alpha})$,
where $\theta_{\alpha}-\alpha\,dz\id_{E_{\alpha}}$
are nilpotent.
Because $\deg(E_{\alpha})\leq 0$
and $\sum_{\alpha} \deg(E_{\alpha})=0$,
we obtain $\deg(E_{\alpha})=0$.
By Theorem \ref{thm;20.9.30.10},
the orthogonal complement of
$E_{\alpha}$ is also Higgs subbundle of $E$,
and hence it is equal to
$\bigoplus_{\beta\neq\alpha}E_{\beta}$.
By an easy induction,
we obtain that
the decomposition $E=\bigoplus E_{\alpha}$
is orthogonal with respect to $h$.
Because we obtain a decomposition of
harmonic bundles
$(E,\delbar_E,\theta,h)
=\bigoplus_{\alpha}
(E_{\alpha},\delbar_{E_{\alpha}},\theta_{\alpha},h_{\alpha})$,
we may assume that
$\theta-\alpha\,dz\,\id_{E}$ is nilpotent
for some complex number $\alpha$
from the beginning.
Moreover, we may assume that $\alpha=0$,
i.e., $\theta$ is nilpotent.
By Lemma \ref{NoSolutionOnParabolicSurface},
we obtain $\theta=0$.
It implies that
the Chern connection of
$(E,\delbar_E,h)$ is flat.
Because the fundamental group of $X$ is abelian,
it is isomorphic to a direct sum of
holomorphic line bundles with a flat Hermitian metric.
\hfill\qed

\subsubsection{Tame harmonic bundles on $\cnum$}

For any positive integer $r$,
we set 
$E(r):=\bigoplus_{i=1}^r\nbigo_{\cnum}\,e_i$.
It is equipped with the trivial Higgs field $0$.
There exits a Hermitian metric $h(r)$
determined by $h(r)(e_i,e_j)=1$ $(i=j)$
or $h(r)(e_i,e_j)=0$ $(i\neq j)$.

\begin{prop}
Any tame harmonic bundle of rank $r$ on $\cnum$
is isomorphic to
$(E(r),0,h(r))$. 
\end{prop}
\pf
Let $(E,\delbar_E,\theta,h)$ be a tame harmonic bundle
on $\cnum$.
We obtain a locally free $\nbigo_{\proj^1}$-module
$\nbigp^h_0\nbige^0$
with the logarithmic Higgs field
$\theta$
as explained in \S\ref{subsection;20.9.29.100}.
We may regard $\theta$
as a section of
$\End(\nbigp^h_0\nbige^0)\otimes K_{\proj^1}(\infty)$.
By considering the product
$K_{\proj^1}(\infty)^{\otimes m}
\otimes
K_{\proj^1}(\infty)^{\otimes \ell}
\simeq
K_{\proj^1}(\infty)^{\otimes (m+\ell)}$
and the product on $\End(\nbigp^h_0\nbige^0)$,
we obtain
a section
$\theta^m$ of
$\End(\nbigp^h_0\nbige^0)\otimes K_{\proj^1}(\infty)^{\otimes\,m}$
for any positive integer $m$.
By taking the trace,
we obtain sections
$\tr(\theta^m)$
of $K_{\proj^1}(\infty)^{\otimes\,m}$.
Because $K_{\proj^1}(\infty)\simeq\nbigo_{\proj^1}(-1)$,
we obtain
$\tr(\theta^m)=0$ for any $m>0$.
It implies that $\theta$ is nilpotent.
By Lemma \ref{NoSolutionOnParabolicSurface},
we obtain $\theta=0$.
It implies that
the Chern connection of
$(E,\delbar_E,h)$ is flat.
Because the fundamental group of $\cnum$ is trivial,
$(E,\delbar_E,h)$ is isomorphic to
$(E(r),h(r))$.
\hfill\qed

\subsubsection{Tame harmonic bundles on $\cnum^{\ast}$}

Let $(a,\alpha)\in\real\times\cnum$.
On $\nbigo_{\cnum^{\ast}}$,
let $h_a$ be the Hermitian metric
determined by
$h_a(1,1)=|z|^{-2a}$.
It is equipped with a Higgs field
$\theta_{\alpha}=\alpha\,dz/z$.
Thus, we obtain a harmonic bundle
$L(a,\alpha)
=\bigl(
 \nbigo_{\cnum^{\ast}},\theta_{\alpha},h_a
 \bigr)$.
 It is easy to see that
 $L(a,\alpha)$ is isomorphic to
 $L(a+n,\alpha)$ for any $n\in\seisuu$.
 
\begin{prop}
For any tame harmonic bundle $(E,\delbar_E,\theta,h)$
of rank $r$
on $\cnum^{\ast}$,
 there exist
 $(a_i,\alpha_i)\in\real\times\cnum$
$(i=1,\ldots,r)$
such that
$(E,\delbar_E,\theta,h) \simeq\bigoplus_{i=1}^r L(a_i,\alpha_i)$. 
\end{prop}
\pf
Let $f$ be the endomorphism of
$E$ determined by $\theta=f\,dz/z$.
We obtain the characteristic polynomial
$\det(t\id_E-f)=\sum_{j=0}^r a_j(z)t^j$.
Because
$f(\nbigp^h_{0,0}\nbige^0)\subset \nbigp^h_{0,0}\nbige^0$,
 $a_j(z)$ are holomorphic at $z=0,\infty$,
and hence constant.
It implies that the eigenvalues of $f$ are constant.
We obtain the decomposition 
$(E,f)=\bigoplus (E_{\alpha},f_{\alpha})$,
where $f_{\alpha}-\alpha\id_{E_{\alpha}}$
are nilpotent.
It induces
the decomposition of the regular filtered Higgs bundles
$(\nbigp^h_{\ast}\nbige^0,\theta)
=\bigoplus (\nbigp^h_{\ast}\nbige^0_{\alpha},\theta_{\alpha})$.
Because
$\deg(\nbigp^h_{\ast}\nbige_{\alpha})\leq 0$
and
$\sum\deg(\nbigp^h_{\ast}\nbige_{\alpha})=0$,
we obtain
$\deg(\nbigp^h_{\ast}\nbige_{\alpha})=0$.
The orthogonal complement
$E_{\alpha}^{\bot}$ is a Higgs subbundle of $E$.
Hence, we obtain
$E_{\alpha}^{\bot}=\bigoplus_{\beta\neq\alpha}E_{\beta}$.
By an easy induction,
we obtain that 
the decomposition $E=\bigoplus E_{\alpha}$
is orthogonal with respect to $h$,
and we obtain the decomposition of harmonic bundle
$(E,\delbar_E,\theta,h)
=\bigoplus_{\alpha\in\cnum}
(E_{\alpha},\delbar_{E_{\alpha}},\theta_{\alpha},h_{\alpha})$.

Hence, we may assume
$\theta-\alpha\,dz/z$ is nilpotent from the beginning.
Moreover,
by considering the tensor product with $L(0,-\alpha)$,
we may assume that $\alpha=0$.
By Lemma \ref{NoSolutionOnParabolicSurface},
we obtain $\theta=0$.
It implies that the Chern connection of
$(E,\delbar_E,h)$ is flat.
Because the fundamental group of $\cnum^{\ast}$ is $\seisuu$,
it is easy to see that
$(E,\delbar_E,h)$ is isomorphic to
$\bigoplus_{i=1}^{\rank(E)} L(a_i,0)$
for some $a_i\in\real$.
\hfill\qed


\section{Estimates for a complete solution}\label{EstimateSection}
In this section, we deduce some quantitative estimates of a complete solution. 

In Corollary \ref{cor;20.9.17.101}, we already see that for a complete solution $h$, the metrics $g(h)_i$ are mutually bounded and $|q|_{g(h)_i}$ is bounded. We will give precise bounds in Theorem \ref{CompletenessRealCurvature}. 

According to Corollary \ref{cor;20.9.18.110}, a complete solution to the system (\ref{eq;20.7.2.1}) must be real. The system (\ref{eq;20.7.2.1}) for a solution $\vecw\in\Toda^{\real}(q,g)$ is reduced to $(w_1,\dots, w_n)(n=[\frac{r}{2}])$ satisfying
\begin{equation}
\label{systemhalf}
\left\{
\begin{array}{l}
\triangle_g w_1=e^{2w_1}|q|^2_g-e^{-w_1+w_2}-\frac{r-1}{4}k_g\\
\triangle_g w_k=e^{-w_{k-1}+w_k}-e^{-w_k+w_{k+1}}-\frac{r+1-2k}{4}k_g,\quad k=2, \cdots, n-1\\
\triangle_g w_n=e^{-w_{n-1}+w_n}-e^{-(2n+2-r)w_n}-\frac{r+1-2n}{4}k_g.
\end{array}
\right.
\end{equation}

Before proving Theorem \ref{CompletenessRealCurvature}, we first prove the following proposition.
\begin{prop}\label{ModelLemma} Suppose that the background metrics $g_1,\cdots, g_m$ are complete and their Gaussian curvature are bounded from below. Suppose $\sigma_1,\cdots, \sigma_m: X\rightarrow \mathbb R$ are $C^2$ functions satisfying the following system
\begin{equation}
\label{ModelSystem}
\left\{
\begin{array}{l}
\triangle_{g_1}\sigma_1\geq (1+a) e^{\sigma_1}-(a+2)+e^{-\sigma_2}\\
\triangle_{g_k}\sigma_k\geq -e^{\sigma_{k-1}+\sigma_k}+3e^{\sigma_k}-3+e^{-\sigma_{k+1}}, \quad k=2, \cdots, m-1\\
\triangle_{g_m}\sigma_m\geq -e^{\sigma_{m-1}+\sigma_{m}}+(2+b)e^{\sigma_{m}}-(1+c),
\end{array}
\right.
\end{equation} where $a, b, c$ are positive constants.

Suppose $d_1,\cdots, d_m$ are positive constants satisfying 
\begin{equation}
\label{ModelSystemConstants}
\left\{
\begin{array}{l}
0= (1+a)d_1-(a+2)+d_2^{-1}\\
0= -d_{k-1}d_k+3d_k-3+d_{k+1}^{-1},\quad k=2, \cdots, m-1\\
0= -d_{m-1}d_m+(2+b)d_m-(1+c).
\end{array}
\right.
\end{equation}
Then either $e^{\sigma_k}<d_k$ for $1\leq k\leq m$ or $e^{\sigma_k}\equiv d_k$ for $1\leq k\leq m$.

If replacing $\sigma_1$ by a nonzero non-negative $C^2$ function $\eta$ and suppose $\eta, \sigma_2,\cdots, \sigma_m: X\rightarrow \mathbb R$ are $C^2$ functions satisfying the following system
\begin{equation}
\label{ModelSystemVariant}
\left\{
\begin{array}{l}
\triangle_{g_1}\eta\geq \eta((1+a)\eta-(a+2)+e^{-\sigma_2})\\
\triangle_{g_2}\sigma_2\geq -\eta e^{\sigma_2}+3e^{\sigma_2}-3+e^{-\sigma_3}\\
\triangle_{g_k}\sigma_k\geq -e^{\sigma_{k-1}+\sigma_k}+3e^{\sigma_k}-3+e^{-\sigma_{k+1}},\quad k=3, \cdots, m-1\\
\triangle_{g_m}\sigma_m\geq -e^{\sigma_{m-1}+\sigma_{m}}+(2+b)e^{\sigma_{m}}-(1+c),
\end{array}
\right.
\end{equation} 
 then we obtain either $\eta< d_1, e^{\sigma_k}<d_k$ for $2\leq k\leq m$ or $\eta\equiv d_1, e^{\sigma_k}\equiv d_k$ for $2\leq k\leq m$. 
\end{prop}

In case $m=1$, the equations (\ref{ModelSystem}), (\ref{ModelSystemConstants}), (\ref{ModelSystemVariant}) are respectively \begin{eqnarray}
&&\triangle_g\sigma\geq (1+a)e^{\sigma}-(a+2), \\
&& 0= (1+a)d-(a+2),\\
&&\triangle_g\eta\geq \eta((1+a)\eta-(a+2)). \end{eqnarray}
\pf
Let $M_i=\sup_{X}e^{\sigma_i}$ for $1\leq i\leq m$. Define the constants $B_i, D_i(1\leq i\leq m+1)$ as follows:
\begin{eqnarray}
B_1=a(1-M_1^{-1}), \quad B_k=2-M_k^{-1}-M_{k-1}(2\leq k\leq m),\quad B_{m+1}=c-bM_m,\\
D_1=a(1-d_1^{-1}), \quad D_k=2-d_k^{-1}-d_{k-1}(2\leq k\leq m),\quad D_{m+1}=c-bd_m.
\end{eqnarray}

\begin{claim}\label{Claim}
(1) $d_1D_1=D_2$; $M_1B_1\leq B_2$ and $M_1\leq \frac{a+2}{a+1}$.\\
(2) $d_kD_k=D_{k+1}$; if $M_{k-1}\leq 2$ where $k$ satisfies $2\leq k\leq m$, then $M_kB_k\leq B_{k+1}$. \\
(3) $d_s-1=\frac{as+1-a}{a}D_{s+1}$; suppose $M_{k-1}\leq 2$ for every $2\leq k\leq s$, then $M_s-1\leq \frac{as+1-a}{a}B_{s+1}$.\\
(4) $M_{k-1}\leq 2$ for every $2\leq k\leq m$.
\end{claim}
\pf (of Claim \ref{Claim})
For Part (1) and (2), the proof will use the Cheng-Yau maximum principle. Note that the background metric $g_1, g_2,\cdots, g_m$ are complete and their Gaussian curvature is bounded from below, the condition of the Cheng-Yau maximum principle on the background metric is satisfied. We only need to check whether the equations satisfy the condition. 

(1) We apply the Cheng-Yau maximum principle to the equation for $\sigma_1$,
\begin{equation}
\triangle_{g_1}\sigma_1\geq(1+a)e^{\sigma_1}-(2+a)+e^{-\sigma_2}\geq (1+a)e^{\sigma_1}-(2+a)+M_2^{-1},
\end{equation} and obtain 
\begin{equation*}
0\geq (1+a)M_1-(2+a)+M_2^{-1}.
\end{equation*}
Note that it is allowed in the equation that $M_k=+\infty$ in which case we write $M_k^{-1}=0$.
So $M_1B_1\leq B_2$ and $M_1\leq \frac{a+2}{a+1}.$ Similarly, $d_1D_1=D_2$.

(2) (i) For some $2\leq k\leq m-1$, suppose $M_{k-1}\leq 2$. The equation for $\sigma_k$ is
\begin{eqnarray}
&&\triangle_{g_k}\sigma_k\geq-e^{\sigma_{k-1}+\sigma_k}+3e^{\sigma_k}-3+e^{-\sigma_{k+1}}\\
\Longrightarrow&&\triangle_{g_k}\sigma_k\geq(3-M_{k-1})e^{\sigma_k}+M_{k+1}^{-1}-3,\label{EquationClaim1}
\end{eqnarray} where the coefficient of $e^{\sigma_k}$ is a positive constant. We apply the Cheng-Yau maximum principle to (\ref{EquationClaim1}) and obtain 
\begin{equation*}
0\geq (3-M_{k-1})M_k+M_{k+1}^{-1}-3,
\end{equation*}
which is equivalent to $M_kB_k\leq B_{k+1}$. Similarly, $d_kD_k=D_{k+1}.$

(ii) Suppose $M_{m-1}\leq 2$, the equation for $\sigma_m$ is
\begin{eqnarray}
&&\triangle_{g_m}\sigma_m\geq-e^{\sigma_{m-1}+\sigma_m}+(2+b)e^{\sigma_m}-(1+c)\\
\Longrightarrow&&\triangle_{g_m}\sigma_m\geq (b+2-M_{m-1})e^{\sigma_m}-(1+c),\label{EquationClaim2}
\end{eqnarray} where the coefficient of $e^{\sigma_m}$ is $(2+b-M_{m-1})$, a positive constant. We apply the Cheng-Yau maximum principle to (\ref{EquationClaim2}) and obtain 
\begin{equation*}
0\geq (b+2-M_{m-1})M_m-(1+c),
\end{equation*}
which is equivalent to $M_mB_m\leq B_{m+1}$. Similarly, $d_mD_m=D_{m+1}.$

(3) By the definitions of $B_i$'s, we have
\begin{eqnarray*}
1-M_s^{-1}&=&B_s+(M_{s-1}-1)\\&=&B_s+M_{s-1}(1-M_{s-1}^{-1})\\
&=&B_s+M_{s-1}B_{s-1}+\cdots+M_{s-1}\cdots M_2B_2+M_{s-1}\cdots M_2(M_1-1)\\
&=&B_s+M_{s-1}B_{s-1}+\cdots+M_{s-1}\cdots M_2B_2+\frac{1}{a}M_{s-1}\cdots M_1B_1.
\end{eqnarray*}
Since $M_1, M_2, \cdots, M_{s-1}\leq 2$, by Part (1) and (2), $M_kB_k\leq B_{k+1}$ for $1\leq k\leq s$ and hence $M_sM_{s-1}\cdots M_kB_k\leq B_{s+1}$ for $1\leq k\leq s$. So
\begin{eqnarray*}
1-M_s^{-1}\leq B_s\cdot (s-1+\frac{1}{a})\leq \frac{as+1-a}{a}M_s^{-1}B_{s+1}\label{EquationBs}
\end{eqnarray*} and thus $M_s-1\leq \frac{as+1-a}{a}B_{s+1}.$

Similarly, we have $d_s-1=\frac{as+1-a}{a}D_{s+1}.$

(4) From Part (1) of Claim \ref{Claim}, $M_1\leq \frac{a+2}{a+1}\leq 2.$ Suppose $s$ be the largest integer such that $M_{k-1}\leq 2$ for $2\leq k\leq s$. Note that $2\leq s\leq m$. 
Applying Part (3) of Claim \ref{Claim}, 
\begin{eqnarray*}&&M_s-1\leq \frac{as+1-a}{a}B_{s+1}=\frac{as+1-a}{a}(2-M_{s+1}^{-1}-M_{s})\leq \frac{as+1-a}{a}(2-M_s)\\\Longrightarrow &&M_s\leq \frac{2as-a+2}{as+1}<2,\end{eqnarray*}which contradicts with the assumption of $s$ unless $s=m$. Therefore, $M_1, M_2, \cdots, M_{m-1}\leq 2$.
\hfill\qed \\
 
Applying Part (3)(4) of Claim \ref{Claim}, 
 \begin{eqnarray*}
 &&M_m-1\leq \frac{am+1-a}{a}B_{m+1}=\frac{am+1-a}{a}(c-bM_m),\\
 &&d_m-1=\frac{am+1-a}{a}D_{m+1}=\frac{am+1-a}{a}(c-bd_m),\end{eqnarray*}
implying that $M_m\leq d_m.$ 
Let $t$ be the smallest integer such that $M_{t+1}\leq d_{t+1}$. Note that $0\leq t\leq m-1$. Applying Part (3)(4) of Claim \ref{Claim},
\begin{eqnarray*}
&&M_t-1\leq \frac{at+1-a}{a}B_{t+1}=\frac{at+1-a}{a}(2-M_{t+1}^{-1}-M_t)\leq\frac{at+1-a}{a}(2-d_{t+1}^{-1}-M_t),\\
&&d_t-1=\frac{at+1-a}{a}D_{t+1}=\frac{at+1-a}{a}(2-d_{t+1}^{-1}-d_t),\end{eqnarray*}
  implying that $M_t\leq d_t$, which contradicts with the assumption of $t$ unless $t=0$. Thus $M_k\leq d_k$ for $1\leq k\leq m$.
 
Applying $M_k\leq d_k$ for $1\leq k\leq m$ into the equation system (\ref{ModelSystem}), 
\begin{equation*}
\left\{
\begin{array}{l}
\triangle_{g_1}\sigma_1\geq (1+a)e^{\sigma_1}-(2+a)+d_2^{-1}\\
\triangle_{g_2}\sigma_k\geq 3e^{\sigma_k}-3-d_{k-1}d_k+d_{k+1}^{-1},\quad 2\leq k\leq m\\
\triangle_{g_m}\sigma_m\geq be^{\sigma_m}-c-d_{m-1}d_m
\end{array}
\right.
\end{equation*}
For each equation, $\log d_k$ is a solution from the definition of $d_k$'s. Applying the strong maximum principle to each equation, we obtain that for every $1\leq k\leq m$, either $e^{\sigma_k}<d_k$ or $e^{\sigma_k}\equiv d_k$. 

If $e^{\sigma_1}\equiv d_1$, we have $0=\triangle_{g_1}\sigma_1\geq(1+a)e^{\sigma_1}-(2+a)+e^{-\sigma_2}=(1+a)\cdot d_1-(2+a)+e^{-\sigma_2}$ implying that $e^{\sigma_2}\equiv d_2.$

If $e^{\sigma_k}\equiv d_k$, we have $0=\triangle_{g_k}\sigma_k\geq-e^{\sigma_{k-1}+\sigma_k}+3e^{\sigma_k}-3+e^{-\sigma_{k+1}}=3d_k-3-e^{\sigma_{k-1}}+e^{-\sigma_{k+1}}$ implying that $e^{\sigma_{k-1}}=d_{k-1}, e^{\sigma_{k+1}}\equiv d_{k+1}.$

If $e^{\sigma_m}\equiv d_m$, we have $0=\triangle_{g_m}\sigma_m\geq-e^{\sigma_{m-1}+\sigma_m}+(2+b)e^{\sigma_m}-(1+c)=-e^{\sigma_{m-1}}d_m+(b+2)d_m-(1+c)$ implying that $e^{\sigma_{m-1}}\equiv d_{m-1}.$

Thus we obtain that either $e^{\sigma_k}<d_k$ for $1\leq k\leq m$ or $e^{\sigma_k}\equiv d_k$ for $1\leq k\leq m$. \\

If we replace $\sigma_1$ by $\eta$, then we apply the Cheng-Yau maximum principle and obtain $\eta$ is bounded from above and $\sup \eta((1+a)\sup \eta-(a+2)+M_2^{-1})\leq 0$. Since $\sup\eta>0$, then $(1+a)\sup \eta-(a+2)+M_2^{-1}\leq 0.$ Let $M_1=\sup \eta$, so we obtain $M_1B_1\leq B_2.$ The rest proof is identical to the above. 
\hfill\qed \\

We don't intend to give a formula of $d_k$'s satisfying (\ref{ModelSystemConstants}) in Proposition \ref{ModelLemma} for general $a, b,c $ since it will be too messy. The following two cases are enough for our later use.
\begin{lem}\label{Constants} Let $d_k, 1\leq k\leq m$, be constants satisfying (\ref{ModelSystemConstants}) in Proposition \ref{ModelLemma}.\\ (1) If $b=c$, then $d_k=1$ for $1\leq k\leq m$; \\(2) If $b=1, c=2$, then $d_k=\frac{(m-k+2)(ma+ka+2-a)}{(m-k+1)(ma+ka+2)}$ for $1\leq k\leq m$. In particular, if $n=[\frac{r}{2}], m=n-1, a=2n+2-r$, then $d_{n-k}=\frac{(k+1)(r-k-1)}{k(r-k)}$ for $1\leq k\leq n-1.$
\end{lem}
\pf
 Applying Part (3)(4) of Claim \ref{Claim},
 \begin{eqnarray*}
 &&d_m-1=\frac{am+1-a}{a}D_{m+1}=\frac{am+1-a}{a}(c-bd_m)\Longrightarrow d_m=\frac{c(am+1-a)+a}{b(am+1-a)+a}.\end{eqnarray*}
 (1) If $b=c$, then $d_m=1$. Let $t$ be the smallest integer such that $d_{t+1}= 1$. Note that $0\leq t\leq m-1$. Applying Part (3) of Claim \ref{Claim},
\begin{eqnarray*}
&&d_t-1=\frac{at+1-a}{a}D_{t+1}=\frac{at+1-a}{a}(2-d_{t+1}^{-1}-d_t)=\frac{at+1-a}{a}(1-d_t),\end{eqnarray*}
  implying that $d_t=1$, which contradicts with the assumption of $t$ unless $t=0$. Thus $d_k=1$ for $1\leq k\leq m$.\\
(2) If $b=1, c=2$, then $d_m=\frac{2ma+2-a}{ma+1}.$ Let $t$ be the smallest integer such that $0\leq t\leq m-1$ and $$d_{t+1}=\frac{(m-(t+1)+2)(ma+(t+1)a+2-a)}{(m-(t+1)+1)(ma+(t+1)a+2)}=\frac{(m-t+1)(ma+ta+2)}{(m-t)(ma+ta+a+2)}.$$ Applying Part (3) of Claim \ref{Claim},
\begin{eqnarray*}
d_t-1&=&\frac{at+1-a}{a}D_{t+1}=\frac{at+1-a}{a}(2-d_{t+1}^{-1}-d_t)\\
&=&\frac{at+1-a}{a}(2-\frac{(m-t)(ma+ta+a+2)}{(m-t+1)(ma+ta+2)}-d_t),\end{eqnarray*}
  implying that $d_t=\frac{(m-t+2)(ma+ta+2-a)}{(m-t+1)(ma+ta+2)}$, which contradicts with the assumption of $t$ unless $t=0$. Thus $d_k=\frac{(m-k+2)(ma+ka+2-a)}{(m-k+1)(ma+ka+2)}$ for $1\leq k\leq m$.
  
Suppose $n=[\frac{r}{2}], m=n-1, a=2n+2-r$. If $r=2n$, then $a=2$ and
\begin{eqnarray*}
d_k=\frac{(m-k+2)(ma+ka+2-a)}{(m-k+1)(ma+ka+2)}=\frac{(n-k+1)(2n+2k-2)}{(n-k)(2n+2k)}=\frac{(n-k+1)(n+k-1)}{(n-k)(n+k)}.
\end{eqnarray*}
If $r=2n+1$, then $a=1$ and
\begin{eqnarray*}
d_k=\frac{(m-k+2)(ma+ka+2-a)}{(m-k+1)(ma+ka+2)}=\frac{(n-k+1)(n+k)}{(n-k)(n+k+1)}.\end{eqnarray*}

In both cases, $d_{n-k}=\frac{(k+1)(r-k-1)}{k(r-k)}$.
\hfill\qed \\

With the preparation of Proposition \ref{ModelLemma} and Lemma \ref{Constants}, we can prove the following estimates for the complete solution. 
\begin{thm}\label{CompletenessRealCurvature}
Suppose $\vecw^{c}=(w_1,\cdots, w_r)\in\Toda(q,g)$ is the complete solution, then one of the following holds:\\
(i)
\begin{eqnarray}
&&\frac{e^{2w_1}|q|_g^2}{e^{-w_1+w_2}}< 1,\quad \frac{(k-1)(r-k+1)}{k(r-k)}<\frac{e^{-w_{k-1}+w_k}}{e^{-w_k+w_{k+1}}}< 1, \quad 2\leq k\leq  n=[\frac{r}{2}],\label{MutualBounded}\\
&&w_k< -\frac{r+1-2k}{r}\log|q|_g,\quad 1\leq k\leq n=[\frac{r}{2}]\label{CompareDifferential}
\end{eqnarray} 
(ii) $w_k=-\frac{r+1-2k}{r}\log|q|_g$ for $1\leq k\leq n$, in which case $q$ has no zeros and $|q|^{\frac{2}{r}}$ defines a complete metric;\\
(iii) $w_k=\log(\frac{(k-1)!}{(r-k)!}2^{r+1-2k})$ for $1\leq k\leq n$, in which case $q\equiv 0$ and $(X, g)$ is a hyperbolic surface.
\end{thm}
\pf
If $q$ has no zeros and $|q|^{\frac{2}{r}}$ defines a complete metric, then the $w_k$'s in Case (ii) is a complete solution. \\
If $q\equiv 0$, then $(X,g)$ has to be a hyperbolic surface for $\Toda(0,g)$ to be nonempty following from Lemma \ref{NoSolutionOnParabolicSurface}. The $w_k$'s in Case (iii) is a complete solution. \\
We only need to show the inequalities in Case (i). \\
I. We first prove the right hand side inequalities of (\ref{MutualBounded}).
Outside zeros of $q$, 
$$\triangle_g(|q|_g^2e^{w})\geq |q|_g^2e^{w}\triangle_g\log(|q|_g^2e^w)=|q|_g^2e^w(\triangle_gw+\frac{r}{2}k_g).$$
Since both sides are continuous on $X$, the above inequality holds on $X$. 
Then we have 
\begin{equation*}
\left\{
\begin{array}{l}
\triangle_g\frac{e^{2w_1}|q|_g^2}{e^{-w_1+w_2}}\geq \frac{e^{2w_1}|q|_g^2}{e^{-w_1+w_2}}(3e^{2w_1}|q|_g^2-4e^{-w_1+w_2}+e^{-w_2+w_3})\\
\triangle_g\log \frac{e^{-w_1+w_2}}{e^{-w_2+w_3}}=-e^{2w_1}|q|_g^2+3e^{-w_1+w_2}-3e^{-w_2+w_3}+e^{-w_3+w_4}\\
\triangle_g\log \frac{e^{-w_{k-1}+w_k}}{e^{-w_k+w_{k+1}}}=-e^{2w_1}|q|_g^2+3e^{-w_{k-1}+w_k}-3e^{-w_k+w_{k+1}}+e^{-w_{k+1}+w_{k+2}},\quad 3\leq k\leq n-1\\
\triangle_g\log \frac{e^{-w_{n-1}+w_n}}{e^{-(2n+2-r)w_n}}=-e^{-w_{n-2}+w_{n-1}}+(2n+4-r)e^{-w_{n-1}+w_n}-(2n+3-r)e^{-(2n+2-r)w_n}
\end{array}
\right.
\end{equation*}
Let $\eta=\frac{e^{2w_1}|q|_g^2}{e^{-w_1+w_2}}, \sigma_2=\log \frac{e^{-w_1+w_2}}{e^{-w_2+w_3}}, \cdots, \sigma_n=\log \frac{e^{-w_{n-1}+w_n}}{e^{-(2n+2-r)w_n}}.$ The above system becomes
\begin{equation}
\label{Beginning}
\left\{
\begin{array}{l}
\triangle_{e^{-w_1+w_2}\cdot g}\eta\geq \eta(3\eta-4+e^{-\sigma_2})\\
\triangle_{e^{-w_2+w_3}\cdot g}\sigma_2=-\eta e^{\sigma_2}+3e^{\sigma_2}-3+e^{-\sigma_3}\\
\triangle_{e^{-w_k+w_{k+1}}\cdot g}\sigma_k=-e^{\sigma_{k-1}+\sigma_k}+3e^{\sigma_k}-3+e^{-\sigma_{k+1}}, \quad 3\leq k\leq n-1\\
\triangle_{e^{-(2n+2-r)w_n}\cdot g}\sigma_n=-e^{\sigma_{n-1}+\sigma_{n}}+(2n+4-r)e^{\sigma_{n}}-(2n+3-r)
\end{array}
\right.
\end{equation}
This is the system (\ref{ModelSystemVariant}) where $m=n, a=2, b=c=2n+2-r.$ By Lemma \ref{Constants}, $d_1=d_2=\cdots=d_n=1$.

Applying Proposition \ref{ModelLemma}, we obtain that either $e^{\sigma_k}<1$ for $1\leq k\leq n$ or $e^{\sigma_k}\equiv 1$ for $1\leq k\leq n$. In the latter case, we obtain that $q$ has no zeros and $2w_1+2\log|q|_g=-w_1+w_2=-w_2+w_3=\cdots=-(2n+2-r)w_n$, implying that $w_k=-\frac{r+1-2k}{r}\log|q|_g$ for $1\leq k\leq n$. And the metric $e^{-w_1+w_2}\cdot g=|q|^{\frac{2}{r}}$ is complete.

II. Next we prove the left hand side inequalities of (\ref{MutualBounded}). We have
\begin{equation*}
\left\{
\begin{array}{l}
\triangle_g\log\frac{e^{-(2n+2-r)w_n}}{e^{-w_{n-1}+w_n}}=e^{-w_{n-2}+w_{n-1}}-(2n+4-r)e^{-w_{n-1}+w_n}+(2n+3-r)w^{-(2n+2-r)w_n}\\
\triangle_g\log \frac{e^{-w_k+w_{k+1}}}{e^{-w_{k-1}+w_k}}=e^{-w_{k-2}+w_{k-1}}-3e^{-w_{k-1}+w_k}+3e^{-w_k+w_{k+1}}-e^{-w_{k+1}+w_{k+2}}, \quad 3\leq k\leq n-1\\
\triangle_g\log \frac{e^{-w_2+w_3}}{e^{-w_1+w_2}}=e^{2w_1}|q|_g^2-3e^{-w_1+w_2}+3e^{-w_2+w_3}-e^{-w_3+w_4}
\end{array}
\right.
\end{equation*}
Let $\eta_1=\log\frac{e^{-(2n+2-r)w_n}}{e^{-w_{n-1}+w_n}}, \eta_2=\log\frac{e^{-w_{n-1}+w_n}}{e^{-w_{n-2}+w_{n-1}}},\cdots, \eta_{n-1}=\log\frac{e^{-w_2+w_3}}{e^{-w_1+w_2}}.$ The above system becomes
\begin{equation}
\label{Later}
\left\{
\begin{array}{l}
\triangle_{e^{-w_{n-1}+w_n}\cdot g}\eta_1= (2n+3-r)e^{\eta_1}-(2n+4-r)+e^{-\eta_2}\\
\triangle_{e^{-w_{n-k}+w_{n-k+1}}\cdot g}\eta_k=-e^{\eta_{k-1}+\eta_k}-3+3e^{\eta_k}-e^{-\eta_{k+1}}, \quad 2\leq k\leq n-2\\
\triangle_{e^{-w_1+w_2}\cdot g}\eta_{n-1}\geq -e^{\eta_{n-2}+\eta_{n-1}}+3e^{\eta_{n-1}}-3
\end{array}
\right.
\end{equation}
This is the system (\ref{ModelSystem}) where $m=n-1, a=2n+2-r, b=1, c=2.$ By Lemma \ref{Constants}, $d_{n-k}=\frac{(k+1)(r-k-1)}{k(r-k)}$ for $1\leq k\leq n-1.$

Applying Proposition \ref{ModelLemma}, we obtain that either $e^{\eta_k}<d_k$ for $1\leq k\leq n$ or $e^{\eta_k}\equiv d_k$ for $1\leq k\leq n$. In the latter case, we obtain that $q\equiv 0$ which cannot happen by assumption.

III. Finally, we prove the inequalities in (\ref{CompareDifferential}). From (\ref{MutualBounded}), we have for each $1\leq k\leq n$,
\begin{eqnarray}
&&(|q|_g^2e^{2w_1})\cdot (e^{-w_1+w_2})^2\cdot \cdots \cdot (e^{-w_{k-1}+w_k})^2<(e^{-w_{k-1}+w_k})^{2k-1}\\
&&<(e^{-w_k+w_{k+1}})^{\frac{(2n+2-r)(2k-1)}{(2n+2-r)(n-k)+1}}\cdots \cdot (e^{-w_{n-1}+w_n})^{\frac{(2n+2-r)(2k-1)}{(2n+2-r)(n-k)+1}}e^{-(2n+2-r)w_n}.
\end{eqnarray}
Hence $|q|_g^2e^{2w_k}<e^{-\frac{(2n+2-r)(2k-1)}{(2n+2-r)(n-k)+1}w_k}$ and thus $w_k<-\frac{r+1-2k}{r}\log|q|_g$.
\hfill\qed \\

Let $f: \widetilde X\rightarrow N:=SL(r,\mathbb C)/SU(r)$ denote the associated equivariant harmonic map for the harmonic bundle $(\hyperk_{X,r}, \theta(q), h^{c})$ where $h^c$ is the harmonic metric corresponding to the complete solution $\vecw^c$. Here, we use the $SL(r,\mathbb C)$-invariant Riemnannian metric on $N$ induced by the Killing form on $\mathfrak{sl}(r,\mathbb C).$
\begin{cor}\label{NegativeCurvature}
Suppose $q$ has at least a zero or $|q|^{2/r}$ does not induce a complete metric on $X$. For each tangent plane of $f(\widetilde X)$, the sectional curvature $K_\sigma^N$ in $N$ satisfies $K_\sigma^N<0$.\\
The curvature $\kappa$ of the pullback metric satisfies
$\kappa<0.$
\end{cor}
\pf
In this proof, $\theta(q)$ and $h^c$ are denoted by $\theta$ and $h$, respectively. By Theorem \ref{CompletenessRealCurvature}, $[\theta,\theta^{\dagger}_h]\neq 0$.
From the curvature formula of $K_\sigma^N$ (see \cite[Proposition 5.3]{LiIntroduction}), 
\[\kappa\leq K_\sigma^N=-\frac{1}{2r}\frac{\bigl|[\theta,\theta^{\dagger}_h]\bigr|_{h,g}^2}
{|\theta|_{h,g}^4-|tr(\theta^2)|_g^2}.\]
Hence  $\kappa\leq K_\sigma^N<0$. 
\hfill\qed

\begin{prop}\label{RealCompletenessBoundedUniqueness}
Suppose $\vecw\in\Toda(q,g)$ is the complete solution. Suppose $(u_1,\cdots, u_n)$ is a supersolution of the system (\ref{systemhalf}) satisfying there exists a constant M such that $w_i-u_i\leq M$ for $1\leq i\leq n$. 

Then either $w_i<u_i$, for $1\leq i\leq n$ or $u_i\equiv w_i$, for $1\leq i\leq n$.
\end{prop}
The definition of a supersolution is in Definition \ref{supersubsolution}.\\
\pf 
Let $\xi_i=w_i-u_i,$ for $1\leq i\leq n$, which are bounded from above. So $\xi_i$'s satisfy
\begin{equation*}
\left\{
\begin{array}{l}
\triangle_{e^{-w_1+w_2}\cdot g}\xi_1\geq(1-e^{-2\xi_1})\frac{e^{2w_1}|q|_g^2}{e^{-w_1+w_2}}-(1-e^{\xi_1-\xi_2}),\\
\triangle_{e^{-w_2+w_3}\cdot g}\xi_2\geq(1-e^{\xi_1-\xi_2})\frac{e^{-w_1+w_2}}{e^{-w_2+w_3}}-(1-e^{\xi_2-\xi_3}),\\
\cdots \\
\triangle_{e^{-(2n+2-r)w_n}\cdot g}\xi_n\geq(1-e^{\xi_{n-1}-\xi_n})\frac{e^{-w_{n-1}+w_n}}{e^{-(2n+2-r)w_n}}-(1-e^{(2n+2-r)\xi_n}).
\end{array}
\right.
\end{equation*}
Denote $M_i=\sup\limits_{X}e^{\xi_i}$ for $1\leq i\leq n$. Since the metric $e^{-w_k+w_{k+1}}\cdot g$ is complete and the curvature is bounded from below, we can apply the Omori-Yau maximum principle to the equation of $\xi_k$. Suppose $s$ is the largest integer such that $M_s=\max\limits_{1\leq i\leq n}M_i$.

If $s=1$, then for every $k$, there exist a sequence of points $p_k\in X$ such that 
$$\triangle_{e^{-w_1+w_2}\cdot g}\xi_1(p_k)\leq \frac{1}{k}, \quad \xi_1(p_k)\geq M_1-\frac{1}{k}.$$
Then at point $p_k$, from the equation of $\xi_1$, 
\begin{eqnarray*}
\frac{1}{k}&\geq& \triangle_{e^{-w_1+w_2}\cdot g}\xi_1(p_k)\geq(1-e^{\xi_r(p_k)-\xi_1(p_k)})\frac{e^{-w_r+w_1}|q|_g^2}{e^{-w_1+w_2}}(p_k)-(1-e^{\xi_1(p_k)-\xi_2(p_k)}).
\end{eqnarray*}
As $k\rightarrow \infty$, we have $0\geq -(1-\frac{M_1}{M_2})$, which contradicts with the definition of $s$.

If $2\leq s\leq n-1$, then for every $k$, there exist a sequence of points $p_k\in X$ such that 
$$\triangle_{e^{-w_s+w_{s+1}}\cdot g}\xi_s(p_k)\leq \frac{1}{k}, \quad \xi_s(p_k)\geq M_s-\frac{1}{k}.$$
Then at point $p_k$, from the equation of $\xi_s$, 
\begin{eqnarray*}
\frac{1}{k}\geq \triangle_{e^{-w_s+w_{s+1}}\cdot g}\xi_s(p_k)&\geq& (1-e^{\xi_{s-1}(p_k)-\xi_s(p_k)})\frac{e^{-w_{s-1}+w_s}}{e^{-w_s+w_{s+1}}}(p_k)-(1-e^{\xi_s(p_k)-\xi_{s+1}(p_k)})\\
&\geq&(1-\frac{M_{s-1}}{M_s-\frac{1}{k}})\delta-(1-\frac{M_s-\frac{1}{k}}{M_{s+1}}).
\end{eqnarray*}
As $k\rightarrow \infty$, we have $0\geq (1-\frac{M_{s-1}}{M_s})\delta-(1-\frac{M_s}{M_{s+1}}),$ which cannot happen since $M_{s+1}<M_s$ and $M_{s-1}\leq M_s$. 

Therefore $s=n$, then for every $k$, there exist a sequence of points $p_k\in X$ such that 
$$\triangle_{e^{-(2n+2-r)w_n}\cdot g}\xi_n(p_k)\leq \frac{1}{k}, \quad \xi_n\geq M_n-\frac{1}{k}.$$
Then at point $p_k$, from the equation of $\xi_n$, 
\begin{eqnarray*}
\frac{1}{k}\geq \triangle_{e^{-(2n+2-r)w_n}\cdot g}\xi_n(p_k)&\geq&(1-e^{\xi_{n-1}(p_k)-\xi_n(p_k)})\frac{e^{-w_{n-1}+w_n}}{e^{-(2n+2-r)w_n}}(p_k)- (1-e^{(2n+2-r)\xi_n(p_k)})\\
&\geq&(1-\frac{M_{n-1}}{M_n-\frac{1}{k}})\delta-(1-(M_n-\frac{1}{k})^{2n+2-r}).
\end{eqnarray*}  
As $k\rightarrow \infty$, we have $0\geq (1-\frac{M_{n-1}}{M_n})\delta-(1-M_n^{2n+2-r}).$ So we have $M_n\leq 1$ and thus $M_i\leq 1$, $1\leq i\leq n$.

Apply $M_i\leq 1$ for $1\leq i\leq n$ to the system and obtain
\begin{equation*}
\left\{
\begin{array}{l}
\triangle_{e^{-w_1+w_2}\cdot g}\xi_1\geq(1-e^{-2\xi_1})e^{2w_1}|q|_g^2-(1-e^{\xi_1}),\\
\triangle_{e^{-w_2+w_3}\cdot g}\xi_2\geq(1-e^{-\xi_2})\frac{e^{-w_1+w_2}}{e^{-w_2+w_3}}-(1-e^{\xi_2}),\\
\cdots \\
\triangle_{e^{-(2n+2-r)w_n}\cdot g}\xi_n\geq(1-e^{-\xi_n})\frac{e^{-w_{n-1}+w_n}}{e^{-(2n+2-r)w_n}}-(1-e^{(2n+2-r)\xi_n}).
\end{array}
\right.
\end{equation*}
Applying the strong maximum principle to each equation, for every $1\leq k\leq n$, either $\xi_k<0$ or $\xi_k\equiv 0$.\\
If $\xi_1\equiv 0,$ then $0=\triangle_{e^{-w_1+w_2}\cdot g}\xi_1\geq(1-e^{-2\xi_1})e^{2w_1}|q|_g^2-(1-e^{\xi_1-\xi_2})=-1+e^{-\xi_2}$, implying that $\xi_2\equiv 0$.\\
If $\xi_k\equiv 0$ for some $2\leq k\leq n-1$, then $0=\triangle_{e^{-w_k+w_{k+1}}\cdot g}\xi_k\geq(1-e^{\xi_{k-1}-\xi_k})\frac{e^{-w_{k-1}+w_k}}{e^{-w_k+w_{k+1}}}-(1-e^{\xi_k-\xi_{k+1}})=(1-e^{\xi_{k-1}})\frac{e^{-w_{k-1}+w_k}}{e^{-w_k+w_{k+1}}}-(1-e^{-\xi_{k+1}})$, implying that $\xi_{k-1}=\xi_{k+1}\equiv 0.$\\
if $\xi_n\equiv 0$, then $0=\triangle_{e^{-(2n+2-r)w_n}\cdot g}\xi_n\geq(1-e^{\xi_{n-1}-\xi_n})\frac{e^{-w_{n-1}+w_n}}{e^{-(2n+2-r)w_n}}-(1-e^{(2n+2-r)\xi_n})=(1-e^{\xi_{n-1}})\frac{e^{-w_{n-1}+w_n}}{e^{-(2n+2-r)w_n}}$, implying that $\xi_{n-1}\equiv 0$.

Hence we obtain either $\xi_k\equiv 0$ for $1\leq k\leq n$ or $\xi_k<0$ for $1\leq k\leq n.$
\hfill\qed

\begin{prop}\label{RealCompletenessUniqueness}
Suppose $\vecw$ is the unique complete solution in $\Toda(q,g)$ and $\vecu\in\Toda^{\real}(q,g)$. Then either $w_i<u_i$, for $1\leq i\leq n=[\frac{r}{2}]$ or $w_i\equiv u_i$, for $1\leq i\leq n=[\frac{r}{2}]$.
\end{prop}
\pf By Corollary \ref{Comparison}, there exists a positive constant $M$ such that 
$$-u_i+u_{i+1}\leq -w_i+w_{i+1}+M, \quad 1\leq i\leq r-1.$$
Since $(w_n-u_n)-(w_{n+1}-u_{n+1})=(2n+2-r)(w_n-u_n)$, then for any integer $i_0$ satisfying $1\leq i_0\leq n$,
\begin{eqnarray*}
w_{i_0}-u_{i_0}&=&\sum\limits_{k=i_0}^{n-1}[(w_k-u_k)-(w_{k+1}-u_{k+1})]+\frac{1}{2n+2-r}[(2n+2-r)(w_n-u_n)]\\
&\leq& (n-i_0)M+\frac{1}{2n+2-r}M\leq nM.\end{eqnarray*}
 Then we obtain the claim of the proposition from Proposition \ref{RealCompletenessBoundedUniqueness}.
\hfill\qed


\section{Existence of a complete real solution} \label{ExistenceSection}
We will show the existence of a complete real solution for the system (\ref{eq;20.7.2.1}), equivalently, a solution to the system (\ref{systemhalf}) satisfying the metrics $e^{-w_1+w_2}\cdot g, \cdots, e^{-w_{n-1}+w_n}\cdot g, e^{-(2n+2-r)w_n}\cdot g$ are complete.

Our main tool is the method of super-subsolution for a system as follows. 

\subsection{Method of super-subsolution for a system}

On a smooth manifold $M$ equipped with a Riemann metric $g$. For an integer $1\leq k\leq n$, let $F_k(x,y_1,\cdots,y_n)$ be a smooth function defined on $M\times \mathbb{R}^n$.  Consider the following system of equations of $(u_1,\cdots, u_n)$
\begin{equation}\label{systemone}
\triangle_g u_k=F_k(x,u_1,\cdots,u_k),\quad 1\leq k\leq n.
\end{equation} Assume that
\begin{equation}\label{DerivativeSign}\frac{\partial F_k}{\partial y_j}\leq 0, \quad\text{ for $j\neq k.$}\end{equation}
In the case $n=1$, $\triangle_g u=f(x, u),$ we do not assume any condition on $f$.

\begin{df} \label{supersubsolution}A vector function $(g_1,\cdots, g_n)$, where $g_k\in C^0(M)\cap W_{loc}^{1,2}(M)$, is called a supersolution (subsolution) of the system (\ref{systemone}) if it satisfies weakly
\begin{equation*}
\triangle_g g_k\leq (\geq) F_k(x,g_1,\cdots,g_n).
\end{equation*}
We denote a supersolution by $\vecw_+$ and a subsolution by $\vecw_-$.
\end{df}

The following existence theorem easily follows from Guest-Lin \cite{GuestLin} with a slight modification. Note here unlike the result in Guest-Lin, we do not need boundary control on the supersolution and subsolution. For two $\mathbb R^n$-valued functions $\vecw$ and $\vecu$, we define $\vecw<\vecu$ (resp. $\vecw\leq \vecu$) if $w_i<u_i$ (resp. $w_i\leq u_i$) for all $1\leq i\leq n$. 

\begin{prop}(Method of super-subsolution)\label{SuperSubExistence}
On a non-compact manifold $M$ equipped with a Riemannian metric $g$, consider the system (\ref{systemone}) which satisfies Assumption (\ref{DerivativeSign}). Suppose $\vecw_-$ and $\vecw_+$ are a subsolution and a supersolution respectively satisfying $\vecw_-<\vecw_+.$ Then there exists a smooth solution $\vecw$ satisfying $\vecw_-\leq \vecw\leq \vecw_+$.
\end{prop}
\pf Let $\{ M_i\}_{i=1}^{\infty}$ be a sequence of compact submanifolds of $M$ with smooth boundary satisfying $ M_1\Subset  M_2\Subset\cdots\Subset  M_n\Subset\cdots$ and $\cup_{i} M_i=M.$  Let $\Phi_i=(\phi_{1,i},\cdots, \phi_{n,i})$ be a smooth vector function over $ M_i$ satisfying $\vecw_-\leq \Phi_i\leq \vecw_+$. The existence of $\Phi_i$ is assured by $\vecw_-<\vecw_+.$  Denote $\vecw_-=(q_1,\cdots, q_n)$ and $\vecw_+=(g_1,\cdots, g_n).$\\

Step 1: We first prove the existence of a solution $\vecw_i$ on $ M_i$ satisfying $\vecw_-\leq \vecw_i\leq \vecw_+$.

On $ M_i$, set $\vecw^{(0)}_i=(w_{1,i}^{(0)},\cdots,w_{n,i}^{(0)}):=(g_1,\cdots,g_n)$, 

On $ M_i$, given a smooth vector function $\vecw^{(l)}_i=(w_{1,i}^{(l)},\cdots,w_{n,i}^{(l)})$, there always exists a smooth vector function $\vecw^{(l+1)}_i=(w_{1,i}^{(l+1)},\cdots,w_{n,i}^{(l+1)})$ (e.g. see \cite{MichaelTaylor} Theorem 1.6 in Chapter 14) satisfying the following linear elliptic Dirichlet problem: for each $1\leq k\leq n$, \begin{equation}\begin{array}{ccc}\label{equation1}
&&\triangle_gw_{k,i}^{(l+1)}-d_{k,i}\cdot w_{k,i}^{(l+1)}=G_k(x,\vecw_i^{(l)}), \quad \text{in $ M_i$}\\
&&w_{k,i}^{(l+1)}=\phi_{k,i},\quad \text{on $\partial  M_i$},
\end{array}\end{equation}
where $d_{k,i}:=\max\Big \{0, \max\limits_{x\in  M_i, q_j(x)\leq y_j\leq g_j(x)} {\frac{\partial F_k}{\partial y_k}}(x,y_1,\cdots,y_n)\Big\}$ is a nonnegative constant and $$G_k(x,y_1,\cdots,y_n):=F_k(x,y_1,\cdots,y_n)-d_{k,i}\cdot y_k.$$

Since (1) $\frac{\partial G_k}{\partial y_k}=\frac{\partial F_k}{\partial y_k}-d_{k,i}\leq 0$ if $q_j\leq y_j\leq g_j$ for all $j$ and (2) $\frac{\partial G_k}{\partial y_j}=\frac{\partial F_k}{\partial y_j}\leq 0$ for $j\neq k$, we obtain $G_k(x,y_1,\cdots,y_n)$ is decreasing in  each $y_i$ if for all $j$, $q_j\leq y_j\leq g_j$. 

Clearly, $\vecw_-\leq \vecw^{(0)}_i= \vecw_+$, then $G_k(\vecw_-)\geq G_k(\vecw^{(0)}_i)= G_k(\vecw_+)$. Since on $\partial M_i$, $q_k\leq w_{k,i}^{(1)}=\phi_{k,i}\leq g_k$, by the weak maximum principle (e.g. see \cite{Trudinger} Theorem 8.1), we obtain $q_k\leq w_{k,i}^{(1)}\leq g_k$ in $ M_i.$ Thus $\vecw_- \leq \vecw^{(1)}_i\leq \vecw^{(0)}_i=\vecw_+.$ Assuming that $\vecw_-\leq \vecw^{(l)}_i\leq \vecw^{(l-1)}_i\leq \vecw_+$, then $G_k(\vecw_-)\geq G_k(\vecw^{(l)}_i)\geq G_k(\vecw^{(l-1)}_i)\geq G_k(\vecw_+)$. Since on $\partial M_i$, $q_k\leq w_{k,i}^{(l+1)}=w_{k,i}^{(l)}=\phi_{k,i}\leq g_k$, by the maximum principle for weak super(-sub)solutions (e.g. see \cite[Theorem 8.1]{Trudinger} ), we conclude $q_k\leq w_{k,i}^{(l+1)}\leq w_{k,i}^{(l)}\leq g_k$ in $ M_i.$ Thus $\vecw_-\leq \vecw^{(l+1)}_i\leq \vecw^{(l)}_i\leq \vecw_+.$ Since $\vecw^{(l)}_i$ is monotone decreasing in $l$ and is bounded, it converges to a function $\widetilde{\vecw}_i$. 

On $ M_{i+3}$, $\widetilde{\vecw}^{(l)}_{i+3}$ are uniformly bounded independent of $l$ and $G_k$ is smooth. In the following context, the constant $C$ varies in different places and does not depend on $l$. Using the interior $L^p$-estimates for linear elliptic equations $(p>\dim M)$ (e.g. see \cite[Theorem 9.11]{Trudinger})
$$||w^{(l)}_{k,i+3}||_{W^{2,p}( M_{i+2})}\leq C(||w^{(l)}_{k,i+3}||_{L^p( M_{i+3})}+||(\triangle_g-d_{k,i+3})w^{(l)}_{k, i+3}||_{L^p( M_{i+3})})\leq C.$$ We conclude the sequence $\{\vecw^{(l)}_{i+3}\}_{l=1}^{\infty}$ is bounded in $W^{2,p}( M_{i+2})$ and hence is bounded in $C^{1,\alpha}( M_{i+1})$ by the Sobolev embedding theorem (e.g. see \cite[Corollary 7.11]{Trudinger}). Note that $||G_k(x, \vecw^{(l-1))}_{i+3})||_{C^{\alpha}( M_{i+1})}\leq C\cdot ||\vecw^{(l-1)}_{i+3}||_{C^{\alpha}( M_{i+2})},$ by the Schauder interior estimates (e.g, see \cite[Theorem 6.2]{Trudinger} ), 
$$||w^{(l)}_{k,i+3}||_{C^{2,\alpha}( M_i)}\leq C(||w^{(l)}_{k,i+3}||_{C^{0}( M_{i+1})}+||(\triangle_g-d_{k,i+3})w^{(l)}_{k,i+3}||_{C^{\alpha}( M_{i+1})}) \leq C.$$ Using the Arzela-Ascoli theorem, the sequence $\{{\vecw}^{(l)}_{i+3}\}_{l=1}^{\infty}$ converges to $\widetilde{\vecw}_{i+3}$ in $C^2( M_i)$. Therefore, the limit $\widetilde{\vecw}_{i+3}$ is a $C^2$ solution of the system (\ref{systemone}) over $ M_i$. By the elliptic regularity and the bootstrap argument, we have that $\widetilde{\vecw}_{i+3}$ is $C^{\infty}$. Moreover, it is clear that $\vecw_-\leq \widetilde{ \vecw}_{i+3}\leq \vecw_+$. Denote $\widetilde{\vecw}_{i+3}$ by $\vecw_i,$ we finish the proof of step 1. \\

Step 2: We prove the existence of solution over $ M.$ On $ M_{i+3}$, the solutions $\vecw_{l+i+3}$ are uniformly bounded independent of $l$ and $F_k$ is smooth. In the following context, the constant $C$ varies in different places and does not depend on $l$. By the boundedness of $\vecw_{l+i+3}$ on $ M_{i+3}$, using the $L^p$-estimates for linear elliptic equations $(p>\dim M)$
$$||w_{k,l+i+3}||_{W^{2,p}( M_{i+2})}\leq C(||w_{k,l+i+3}||_{L^p( M_{i+3})}+||\triangle_g w_{k,l+i+3}||_{L^p( M_{i+3})})\leq C,$$ 
we conclude the sequence $\{\vecw_{l+i+3}\}_{l=1}^{\infty}$ is bounded in $W^{2,p}( M_{i+2})$ and hence is bounded in $C^{1,\alpha}( M_{i+1})$ by the Sobolev embedding theorem. Note that $||F_k(x, \vecw_{l+i+3})||_{C^{\alpha}( M_{i+1})}\leq C\cdot ||\vecw_{l+i+3}||_{C^{\alpha}( M_{i+2})}.$ By the Schauder interior estimates (e.g, see \cite[Theorem 6.2]{Trudinger}), 
$$||\vecw_{l+i+3}||_{C^{2,\alpha}( M_i)}\leq C(||\vecw_{l+i+3}||_{C^{\alpha}( M_{i+1})}+||\triangle_g \vecw_{l+i+3}||_{C^{\alpha}( M_{i+1}}) \leq C=C(i).$$  Using the Arzela-Ascoli theorem and a diagonal trick, we can find a subsequence $\{{\vecw}_{l'+i+3}\}_{l'=1}^{\infty}$ which converges to ${\vecw}$ in $C_{loc}^2(  M)$, so that $ {\vecw}$ is a $C^2$ solution of the system. It follows from the elliptic regularity and the bootstrap argument that $ {\vecw}$ is indeed a $C^{\infty}$ solution of the system over $ M$.
\hfill\qed

\begin{rem}
We may compare the method of super-subsolution developed in Proposition \ref{SuperSubExistence} with Theorem 9 in Wan \cite{WAN} for a scalar equation $\triangle_gu=f(x,u)$. The difference here is that (1) we do not impose the condition $\frac{\partial f}{\partial u}\geq 0$; (2) the non-compact manifold $(M,g)$ is not necessarily complete. 
\hfill\qed 
\end{rem}

\begin{lem}\label{Properties}
(1) Suppose $(\psi_1,\cdots,\psi_n)$ and $(\eta_1,\cdots,\eta_n)$ are two supersolutions. Define $g_k=\min\{\psi_k,\eta_k\}$, $k=1,\cdots,n$, then $(g_1,\cdots,g_n)$ is again a supersolution.\\
(2) Suppose $(\xi_1,\cdots,\xi_n)$ and $(\phi_1,\cdots,\phi_n)$ are two subsolutions. Define $q_k=\max\{\xi_k,\phi_k\}$, $k=1,\cdots,n$, then $(q_1,\cdots,q_n)$ is again a subsolution.
\end{lem}
\pf
(1) Suppose $(\psi_1,\cdots,\psi_n)$ and $(\eta_1,\cdots,\eta_n)$ are two supersolutions. 
Since $\frac{\partial F_k}{\partial u_j}\leq 0, j\neq k$ and $g_k\leq \psi_k, \eta_k$, then for $g_k$ we have \begin{equation*}\triangle_g\psi_k\leq F_k(x,\psi_1,\psi_2,\cdots,\psi_n)\leq F_k(x,g_1,\cdots,g_{k-1},\psi_k,g_{k+1},\cdots,g_n),\end{equation*}
\begin{equation*}\triangle_g\eta_k\leq F_k(x,\eta_1,\eta_2,\cdots,\eta_n)\leq F_k(x,g_1,\cdots,g_{k-1},\eta_k,g_{k+1},\cdots,g_n).\end{equation*}
Since $g_k=\min \{\psi_k, \eta_k\}$ and using the fact that the minimum of two supersolutions is still a supersolution (see \cite{DancerSweers} for example), one gets 
$$\triangle_g g_k\leq F_k(x,g_1,\cdots,g_n).$$
Hence $(g_1,\cdots,g_n)$ is a supersolution.

(2) Suppose $(\xi_1,\cdots,\xi_n)$ and $(\phi_1,\cdots,\phi_n)$ are two subsolutions. Then $(-\xi_1,\cdots,-\xi_n)$ and $(-\phi_1,\cdots,-\phi_n)$ are two supersolutions to the system
\begin{equation}\label{OppositeSystem}
\triangle_g u_k=G_k(u_1,\cdots, u_n):=-F_k(x,-u_1,\cdots,-u_n).
\end{equation}
Clearly $G_k$ satisfies the Assumption (\ref{DerivativeSign}). By Part (1), $(-q_1,\cdots, -q_n)$ is a supersolution to the system (\ref{OppositeSystem}) and thus it is a subsolution to the system (\ref{systemone}).
\hfill\qed 

\subsection{Existence of a complete real solution}

In this subsection, we will show the existence of a complete real solution to the system (\ref{eq;20.7.2.1}) for a Riemann surface $X$ and a holomorphic $r$-differential $q$. The plan of the proof goes as follows. 

1. We first show the existence in Proposition \ref{BoundedExistenceBackground} when $X$ is the unit disk $\mathbb D$ with an extra condition that $|q(z)|(1-|z|^2)^r$ is bounded. With this extra condition, it is easy to construct a supersolution and a subsolutiom and we can apply the method of super-subsolution.

2. Next, we remove the condition that $|q(z)|(1-|z|^2)^r$ is bounded and show the existence in Proposition \ref{Existence} for any $q$ on the unit disk $\mathbb D$. We start with an exhaustion of the unit disk by a sequence of disks. Since on each smaller disk, the boundedness condition on $q$ is automatically satisfied and the existence follows from Step 1. Then the main difficulty is to show the sequence of solutions  converges in $C^2_{loc}(\mathbb D)$. 

3. We then show the existence in Proposition \ref{HyperbolicRiemannSurface} when $X$ is a hyperbolic Riemann surface, that is, it has a universal cover as the unit disk $\mathbb D$. By lifting $(q, g)$ to the universal cover $\mathbb D$, we obtain a complete solution for the lifted pair from Step 2. The solution will descend to a solution on $X$ because of the uniqueness property of a complete solution. 

4. We show the existence in Proposition \ref{ComplexPlaneExistence} when $X$ is the complex plane $\mathbb C$. This is again done by constructing a supersolution and a subsolution. Unlike Step 1, the supersolution and subsolution are not easy to be constructed directly. We cover the complex plane with a disk and the complement of a smaller disk. On each piece, we obtain a complete solution by Step 3 since both of them are hyperbolic Riemann surfaces. The supersolution and subsolution arise from a gluing of these two solutions.

5. Using a similar method as Step 3, we then show the existence in Proposition \ref{ParabolicRiemannSurface} when $X$ is a parabolic Riemann surface, that is, it has a universal cover as the complex plane $\mathbb C$.

In this way, we finish proving the existence of a complete real solution for any Riemann surface $X$ and a holomorphic $r$-differential $q$ over $X.$ In the end of the subsection, we also prove the existence of a real solution in Proposition \ref{FiniteManyZerosExistence} if $q$ has finitely many zeros which is asymptotic to the behavior of $|q|$. This will lead to a discussion of the non-uniqueness of solutions if we remove the completeness condition.\\

When $X=\mathbb C$ or any subdomain of $\mathbb C$, by choosing the background K\"ahler metric as $dz\otimes d\bar z$, the system (\ref{systemhalf}) is just looking for a solution $(\widetilde w_1,\cdots, \widetilde w_n)$ satisfying
\begin{equation}
\label{systemhalfbackground}
\left\{
\begin{array}{l}
\triangle\widetilde w_1=e^{2\widetilde w_1}|q|^2-e^{-\widetilde w_1+\widetilde w_2}\\
\triangle \widetilde w_2=e^{-\widetilde w_1+\widetilde w_2}-e^{-\widetilde w_2+\widetilde w_3}\\
\cdots\\
\triangle \widetilde w_n=e^{-\widetilde w_{n-1}+\widetilde w_n}-e^{-(2n+2-r)\widetilde w_n},
\end{array}
\right.
\end{equation} where $\triangle=\partial_z\partial_{\bar z}.$
We denote $(\widetilde w_1,\cdots, \widetilde w_n)$ by $\widetilde{\vecw}$. From now on, if we use $\widetilde \vecw$ instead of $\vecw$, that means we are dealing with the background metric $dz\otimes d\bar z$. 

Consider a domain $\Omega$ equipped with a hyperbolic K\"ahler metric $g=g_0dz\otimes d\bar z$, that is, $\triangle\log g_0=\frac{1}{2}g_0$ for $\triangle=\partial_z\partial_{\bar z}$. If $q=0$, the system (\ref{systemhalfbackground}) admits a unique complete solution $(\widetilde w_1,\cdots, \widetilde w_n)$ where $$\widetilde w_l=\log\Big(\frac{(l-1)!}{(r-l)!}\cdot(\frac{g_0}{4})^{-\frac{r+1-2l}{2}}\Big ),\quad 1\leq l\leq n.$$ We denote the solution by $\widetilde{\vecw}_{base}$. In particular, on the unit disk $\mathbb D$, $g_0=\frac{4}{(1-|z|^2)^2}$. Each entry of $\widetilde{\vecw}_{base}$ goes to $-\infty$ as $z\rightarrow \partial \Omega$. Clearly, $\widetilde{\vecw}_{base}$ is a supersolution for general $q$. 

On a general domain $\Omega$, if $q\neq 0$, define the vector function
$$\widetilde{\vecw}_{q}=(-\frac{r-1}{r}\log|q|,-\frac{r-3}{r}\log |q|,\cdots,-\frac{r+1-2n}{r}\log|q|),$$ which is a solution to the system (\ref{systemhalfbackground}) except at zeros of $q$. Moreover, each entry of $\widetilde{\vecw}_{q}$ goes to $+\infty$ at each zero of $q.$

\begin{rem}
(i) For a domain $\Omega\subset \mathbb C$,  we may equip two natural K\"ahler metrics: $dz\otimes d\bar z$ and $g_0dz\otimes d\bar z$. For example, on a disk $\mathbb D$, $g_0=\frac{4}{(1-|z|^2)^2}$. Let $\vecw=(w_1,\cdots, w_n)$ where $e^{w_l}=e^{\widetilde w_l}
\cdot g_0^{\frac{r+1-2l}{2}} (1\leq l\leq n)$. Then $\widetilde\vecw$ is a solution to (\ref{systemhalfbackground}) on $(\Omega, dz\otimes d\bar z)$ if and only if $\vecw$ is a solution to (\ref{systemhalf}) on $(\Omega, g_0dz\otimes d\bar z)$. \\
(ii) Suppose $\Omega$ is a hyperbolic domain. Denote $\vecw_{base}$ where $w_l=\log(\frac{(l-1)!}{(r-l)!}2^{(r+1-2l)})$ for $1\leq l\leq n$, which is a solution to (\ref{systemhalf}) for $q=0$ case. Then $\vecw_{base}$ is a solution to (\ref{systemhalf}) which corresponds to the solution $\widetilde\vecw_{base}$ to (\ref{systemhalfbackground}) for $q=0$ case.\\
(iii) Define $\vecw_{q}$ where $w_l=-\frac{r+1-2l}{r}\log |q|_g$ for $1\leq l\leq n$, which is a solution to (\ref{systemhalf}) except at zeros of $q$. Then $\vecw_{q}$ is a solution to (\ref{systemhalf}) which corresponds to the solution $\widetilde\vecw_{q}$ to (\ref{systemhalfbackground}) except at zeros of $q$.
\hfill\qed 
\end{rem}

\begin{prop}\label{BoundedExistenceBackground}
On the unit disk $\mathbb D$, suppose $q(z)$ is a  holomorphic $r$-differential on $\mathbb D$ satisfying $|q(z)|(1-|z|^2)^r$ is bounded, then there exists a smooth solution $\widetilde {\vecw}=(\widetilde w_1,\cdots, \widetilde w_n)$ of the system (\ref{systemhalfbackground}) satisfying 
\[\widetilde{\vecw}_{base}-c\leq \widetilde{\vecw}\leq \widetilde{\vecw}_{base}(\leq 0),\]
where $c$ is a positive constant only depending on $\sup_{\mathbb D} |q(z)|(1-|z|^2)^r$ and $r$, and $\vecu-c:=(u_1-c, \cdots, u_n-c)$.

Moreover, $\widetilde{\vecw}$ is a complete solution satisfying 
\begin{eqnarray}
&&e^{-\widetilde w_l+\widetilde w_{l+1}}\geq l(r-l)\frac{1}{(1-|z|^2)^2}\quad(1\leq l\leq n-1), \quad e^{-(2n+2-r)\widetilde w_n}\geq n(r-n)\frac{1}{(1-|z|^2)^2},\label{BoundedInequalities}\\
&&e^{-\widetilde w_l}\geq \frac{(r-l)!}{(l-1)!}(\frac{1}{(1-|z|^2)^2})^{(r+1-2l)/2}\quad(1\leq l\leq n).\label{BoundedInequalities2}
\end{eqnarray}
\end{prop}
\pf Let $g_0=\frac{4}{(1-|z|^2)^2}$, which satisfies $\triangle \log g_0=\frac{1}{2}g_0$ for $\triangle=\partial_z\partial_{\bar z}$. Then $g=g_0dz\otimes d\bar z$ is the unique complete hyperbolic K\"ahler metric on $\mathbb D$. By assumption $|q|_g=|q(z)|(\frac{1-|z|^2}{2})^r$ is bounded.

Clearly, $\widetilde\vecw_+=\widetilde{\vecw}_{base}$ is a supersolution. Next we construct a subsolution. Define an $n$-tuple of negative constants $(a_1,\cdots, a_n)$ satisfying 
\begin{equation*}
{a_{k+1}-a_k}=\log (kE)(1\leq k\leq n-1),\quad a_n=-\log (nE),
\end{equation*}
where $E$ is a positive constant to be determined later. Define $\xi_k=a_k+\log (\frac{g_0}{4})^{-\frac{r+1-2l}{2}}.$ Let $M=\sup_{\mathbb D} |q|_{g}^2$. For $E$ large enough, $(\xi_1,\cdots,\xi_n) $ satisfies 
\begin{equation*}
\left\{
\begin{array}{l}
-\frac{r-1}{2}\triangle\log g_0=\triangle\xi_1\geq M\cdot e^{2\xi_1}-e^{\xi_2-\xi_1}\geq (e^{2a_1}|q|_{g}^2-e^{a_2-a_1})\cdot \frac{g_0}{4}\\
-\frac{r-3}{2}\triangle\log g_0=\triangle\xi_2\geq e^{\xi_2-\xi_1}-e^{\xi_3-\xi_2}=(e^{a_2-a_1}-e^{a_3-a_2})\cdot \frac{g_0}{4}\\
\cdots\\
-\frac{r+1-2n}{2}\triangle\log g_0=\triangle\xi_n\geq e^{\xi_n-\xi_{n-1}}-e^{-(2n+2-r)\xi_n}=(e^{a_n-a_{n-1}}-e^{-(2n+2-r)a_n})\cdot \frac{g_0}{4}.
\end{array}
\right.
\end{equation*}
Thus $\widetilde\vecw_-=(\xi_1,\cdots,\xi_n)$ is a subsolution.

As long as $E$ is large enough, we have $\widetilde\vecw_-< \widetilde\vecw_+$. Applying Proposition \ref{SuperSubExistence}, there is a $C^{\infty}$ solution $\widetilde{\vecw}$ satisfying $\widetilde\vecw_-\leq \widetilde{\vecw}\leq \widetilde\vecw_+$. 

For $1\leq k\leq n-1$, $e^{-\widetilde w_k+\widetilde w_{k+1}}|dz|^2\geq e^{-a_{k+1}}\cdot (\frac{(r-k)!}{(k-1)!})\cdot \frac{g_0}{4}|dz|^2$ which is a complete metric; $e^{-(2n+2-r) \widetilde w_n}|dz|^2\geq (\frac{(r-n)!}{(n-1)!})^{2n+2-r}\cdot \frac{g_0}{4}|dz|^2$ which is a complete metric. Hence $\widetilde{\vecw}$ is a complete solution and by Proposition \ref{CompletenessRealCurvature}, it follows that the solution $\widetilde{\vecw}$ satisfies
\begin{eqnarray*}
\frac{(k-1)(r-k+1)}{k(r-k)}<\frac{e^{-\widetilde w_{k-1}+\widetilde w_k}}{e^{-\widetilde w_k+\widetilde w_{k+1}}}(2\leq k\leq  n-1),\quad
\frac{(n-1)(r-n+1)}{n(r-n)}< \frac{e^{-\widetilde w_{n-1}+\widetilde w_{n}}}{e^{-(2n+2-r)\widetilde w_n}}.
\end{eqnarray*} 
Since $e^{-(2n+2-r)\widetilde w_n}\geq (\frac{(r-n)!}{(n-1)!})^{2n+2-r}\frac{1}{(1-|z|^2)^2}=n(r-n)\frac{1}{(1-|z|^2)^2}$, together with the above inequalities, we obtain the required estimates in (\ref{BoundedInequalities})(\ref{BoundedInequalities2}).
\hfill\qed \\

Now we show the existence of a complete solution without the condition that $|q(z)|(1-|z|^2)^r$ is bounded. We will use the method developed in Wan \cite{WAN} by using the exhaustion of the unit disk for a scalar equation. Some steps there are simplified. 
\begin{prop}\label{Existence}
Let $q$ be a holomorphic $r$-differential on the disk $\mathbb D$, then the system (\ref{systemhalfbackground}) admits a complete $C^{\infty}$-solution $\widetilde{\vecw}=(\widetilde w_1,\cdots, \widetilde w_n)$. Moreover, the solution satisfies Inequality (\ref{BoundedInequalities})(\ref{BoundedInequalities2}).
\end{prop}
\pf
For any positive integer $k$, let $D_k=\{z\in \mathbb D, |z|<R_k\},$ where
$R_k=\frac{e^k-1}{e^k+1}<1$. Then we have $D_1\subset \cdots \subset D_k\subset D_{k+1}\subset \cdots\subset \mathbb D$, and $\mathbb D=\cup_{k=1}^{\infty}D_k.$ On $D_k$, it has the Poincar\'e metric $\frac{4R_k^2}{(R_k^2-|z|^2)^2}dz\otimes d\bar z.$ On each $D_k$, $|q|$ is bounded and hence $|q|(R_k^2-|z|^2)^r$ is bounded. By Proposition \ref{BoundedExistenceBackground}, for each $k$, there is a complete solution $\widetilde {\vecw}^{(k)}$ to the system (\ref{systemhalfbackground}) on $D_k$. 

Since $\widetilde {\vecw}^{(k+1)}$ restricts to a bounded solution on $D_k$ and $\widetilde {\vecw}^{(k)}$ goes to $-\infty$ as approaching to $\partial D_k$, by Proposition \ref{RealCompletenessBoundedUniqueness}, we have $\widetilde {\vecw}^{(k+1)}\geq \widetilde {\vecw}^{(k)}$ on $D_k$. Since $\widetilde{\vecw}_{base}$ restricts to a supersolution on each $D_k$ and $\widetilde {\vecw}^{(k)}$ goes to $-\infty$ as approaching to $\partial D_k$, by Proposition \ref{RealCompletenessBoundedUniqueness}, the complete solution $\widetilde{\vecw}^{(k)}$ satisfies $\widetilde{\vecw}^{(k)}\leq \vecw_{base}.$ We conclude the sequence $\{\widetilde {\vecw}^{(k)}\}_k$ is monotone increasing and is bounded from above. Hence the sequence $\widetilde {\vecw}^{(k)}$ converges pointwise to a vector function $\widetilde {\vecw}$. 

On $\overline D_{k+2}$, the sequence $\{\widetilde{\vecw}^{(k+j+3)}\}_{j=1}^{\infty}$ is uniformly bounded independent of $j$. In the following context, the constant $C=C(k)$ varies in different places and does not depend on $j$. Using the interior $L^p$-estimates for linear elliptic equations $(p>2)$ (e.g. see \cite{Trudinger} Theorem 9.11)
$$||\widetilde{\vecw}^{(k+j+3)}||_{W^{2,p}(\overline D_{k+2})}\leq C(||\widetilde{\vecw}^{(k+j+3)}||_{L^p(\overline D_{k+3})}+||\triangle \widetilde{\vecw}^{(k+j+3)}||_{L^p(\overline D_{k+3})})\leq C,$$ we conclude the sequence $\{\widetilde {\vecw}^{(k+j+3)}\}_{j=1}^{\infty}$ is bounded in $W^{2,p}(\overline D_{k+2})$ and is also bounded in $C^{1,\alpha}(\overline D_{k+1})$ by the Sobolev embedding theorem (e.g. see \cite{Trudinger} Corollary 7.11). Since $|F_k(x, \widetilde{\vecw}^{(k+j+3)})|_{C^{\alpha}(\overline D_{k+1})}\leq C\cdot ||\widetilde{\vecw}^{(k+j+3)}||_{C^{\alpha}(\overline D_{k+1})},$ by the Schauder interior estimates (e.g, see \cite{Trudinger} Theorem 6.2), 
$$||\widetilde{\vecw}^{(k+j+3)}||_{C^{2,\alpha}(\overline  D_k)}\leq C(||\widetilde{\vecw}^{(k+j+3)}||_{C^0(\overline D_{k+1})}+||\triangle \widetilde{\vecw}^{(k+j+3)}||_{C^{\alpha}(\overline D_{k+1}}) \leq C=C(k).$$ Using the Arzela-Ascoli theorem and a diagonal trick, we can find a subsequence $\{\widetilde {\vecw}^{(k+j'+3)}\}_{j'=1}^{\infty}$ which converges to $\widetilde {\vecw}$ in $C_{loc}^2(\mathbb D)$, so that $\widetilde {\vecw}$ is a $C^2$ solution of the system. By the elliptic regularity and the bootstrap argument, we have that $\widetilde {\vecw}$ is indeed a $C^{\infty}$ solution of the system.

Finally, by Proposition \ref{BoundedExistenceBackground}, for each $k$, we have 
\begin{eqnarray*}
e^{-\widetilde w_l^{(k)}+\widetilde w_{l+1}^{(k)}}\geq l(r-l)\frac{R_k^2}{(R_k^2-|z|^2)^2}\quad(1\leq l\leq n-1),\quad
e^{-(2n+2-r)\widetilde w_n^{(k)}}\geq n(r-n)\frac{R_k^2}{(R_k^2-|z|^2)^2}.
\end{eqnarray*}
It follows that \begin{eqnarray*}
e^{-\widetilde w_l+\widetilde w_{l+1}}\geq l(r-l)\frac{1}{(1-|z|^2)^2}\quad(1\leq l\leq n-1),\quad
e^{-(2n+2-r)\widetilde w_n}\geq n(r-n)\frac{1}{(1-|z|^2)^2}.
\end{eqnarray*} Thus the solution $\widetilde {\vecw}$ is complete. The estimates follow from the above equalities.
\hfill\qed 

\begin{prop}\label{HyperbolicRiemannSurface}
Let $(X, g)$ be a complete hyperbolic surface. Let $q$ be a holomorphic $r$-differential on $X$, then the system (\ref{systemhalf}) admits a complete solution $\vecw$. 

Moreover, 
\begin{eqnarray}
&&e^{- w_l+ w_{l+1}}\geq \frac{l(r-l)}{4}\quad(1\leq l\leq n-1), \quad e^{-(2n+2-r)w_n}\geq \frac{n(r-n)}{4},\label{BoundedInequalitiesGeneral}\\
&&w_l\leq\log(\frac{(l-1)!}{(r-l)!}2^{(r+1-2l)})\quad(1\leq l\leq n).\label{BoundedInequalitiesGeneral2}
\end{eqnarray}

In particular, if $q$ is bounded with respect to $g$, then the system (\ref{systemhalf}) admits a bounded solution $\vecw$.
\end{prop}
\pf
Let $X$ be covered by $\mathbb D$ under the map $p:\mathbb D\rightarrow X$, with the covering transformation group of $X$ be $\Gamma<Aut(\mathbb D)=PSL(2,\mathbb R)$, i.e. $X=\mathbb D/\Gamma$. Lift $q, g$ to $\hat q, \hat g$ on $\mathbb D$, which are invariant under $\Gamma$. By Proposition \ref{Existence}, there exists a complete solution $\widetilde\vecw\in \Toda(\hat q, dz\otimes d\bar z)$, equivalently, there exists a complete solution $\hat{\vecw}\in \Toda(\hat q, \hat g)$. By Proposition \ref{Pullback}, $\gamma^*(\hat\vecw)=\hat\vecw,$ for $\gamma\in \Gamma$. Hence $\hat {\vecw}$ descends to a solution $\vecw$ on $\mathbb D/\Gamma=X$ of the system which is still complete. The estimates directly follow from (\ref{BoundedInequalities})(\ref{BoundedInequalities2}). 

If $q$ is bounded with respect to $g$, by Proposition \ref{BoundedExistenceBackground}, there exists a bounded solution $\hat {\vecw}\in \Toda(\hat q, \hat g)$ which is complete. By Proposition \ref{Pullback}, $\gamma^*(\hat\vecw)=\hat\vecw,$ for $\gamma\in \Gamma$. Hence $\hat {\vecw}$ descends to a solution $\vecw$ on $\mathbb D/\Gamma=X$ of the system which is still bounded. 
\hfill\qed

\begin{cor}\label{BoundedBounded}
Let $(X, g)$ be a complete hyperbolic surface. If an $r$-differential $q$ on $X$ is bounded with respect to $g$, there uniquely exists a bounded solution in $\Toda(q,g)$. Moreover, it is real. Conversely, if there exists a bounded solution in $\Toda(q,g)$, then $q$ is bounded with respect to $g$. 

In fact, the bounded solution is the complete solution. 
\end{cor}
\pf 
If an $r$-differential $q$ on $X$ is bounded with respect to $g$, by Proposition \ref{HyperbolicRiemannSurface} and Corollary \ref{BoundedSolutionUnique}, there exists a unique bounded solution. Conversely, if there exists a bounded solution in $\Toda(q,g)$, by Corollary \ref{BoundedSolutionImplyBoundedDifferential}, $q$ is bounded with respect to $g$.

Obviously, the bounded solution is complete.
\hfill\qed\\

Now we move to the case of complex plane. We will use the method developed in Au-Wan \cite{AuWAN} to construct an appropriate supersolution and subsolution.
\begin{prop}\label{ComplexPlaneExistence}
Let $q$ be a nonzero holomorphic $r$-differential on $\mathbb C$, then the system (\ref{systemhalfbackground}) admits a complete solution. 
\end{prop}
\pf
We will use the method of super-subsolution for system to show the existence by constructing a supersolution and a subsolution in the weak sense. 

Without loss of generality, we assume $q(0)\neq 0$. Choose $R_1$ such that the zeros of $q$ are outside the ball $B_{R_1}$. Choose $R_2<R_1$. Let $\widetilde\vecu$ be the unique complete solution to (\ref{systemhalfbackground}) on $B_{R_1}$ and $\widetilde\vecv$ be the unique complete solution to (\ref{systemhalfbackground}) in the complement of the closed disk $\overline{B}_{R_2}$ as described in Proposition \ref{HyperbolicRiemannSurface}.

We first construct a subsolution $\widetilde\vecw_{-}$ as follows:\\
on $\{|z|<R_2\}$, let $\widetilde\vecw_{-}=\widetilde\vecu$;\\
on $\{R_2\leq |z|\leq R_1\}$, let $\widetilde\vecw_{-}=\max\{\widetilde\vecu,\widetilde\vecv\}$;\\
on $\{|z|>R_1\}$, let $\widetilde\vecw_{-}=\widetilde\vecv$.\\
By Proposition \ref{HyperbolicRiemannSurface}, $e^{-\frac{2}{r+1-2l}\widetilde u_l}|dz|^2$ dominates the hyperbolic metric on $B_{R_1}$ up to a constant and hence blows up at boundary of $B_{R_1}$. Then at the neighborhood of $\partial B_{R_1}$, $\widetilde\vecw_{-}=\widetilde\vecv$. Similarly, by Proposition \ref{HyperbolicRiemannSurface}, $e^{-\frac{2}{r+1-2l}\widetilde v_l}|dz|^2$ dominates the hyperbolic metric on the complement of $B_{R_2}$ up to a constant and hence blows up at boundary of $B_{R_2}$. Then at the neighborhood of $\partial B_{R_2}$, $\widetilde\vecw_{-}=\widetilde\vecu$. Hence $\widetilde\vecw_-$ is continuous. Since both $\widetilde\vecu$ and $\widetilde\vecv$ are solutions on the annulus $\{R_2\leq |z|\leq R_1\}$, by Proposition \ref{Properties}, the maximum of $\widetilde\vecu$ and $\widetilde\vecv$ is a subsolution in the weak sense. Therefore $\widetilde\vecw_-\in C^0(\mathbb C)\cap W^{1,2}_{loc}(\mathbb C)$ is a subsolution.

Secondly, we construct a supersolution $\widetilde\vecw_+$. Choose $R_1', R_2'$ such that $R_2<R_2'<R_1'<R_1$. For a constant $c$ and a solution $\widetilde\vecv$ on a domain $U\subset \mathbb C$, $\widetilde\vecv+c=(\widetilde v_1+c, \widetilde v_2+c,\cdots \widetilde v_n+c)$ is a supersolution (subsolution) on $U$ for $c\geq 0$ ($c\leq 0$). In fact, for $c\geq 0$ ($c\leq 0$),
\begin{equation*}
\left\{
\begin{array}{l}
\triangle (\widetilde v_1+c)=\triangle \widetilde v_1\leq (\geq) e^{2(\widetilde v_1+c)}|q|^2-e^{-(\widetilde v_1+c)+(\widetilde v_2+c)}\\
\triangle (\widetilde v_2+c)=\triangle \widetilde v_2=e^{-(\widetilde v_1+c)+(\widetilde v_2+c)}-e^{-(\widetilde v_2+c)+(\widetilde v_3+c)}\\
\quad\quad\quad\cdots\\
\triangle (\widetilde v_n+c)=\triangle \widetilde v_n\leq (\geq) e^{-(\widetilde v_{n-1}+c)+(\widetilde v_n+c)}-e^{-2(\widetilde v_n+c)}.
\end{array}
\right.
\end{equation*}
On the annulus $\{R_2'\leq |z|\leq R_1'\}$, $\widetilde v_l$ has a minimum and $-\frac{r+1-2l}{r}\log |q|$ has a maximum since all the zeros of $q$ are outside $B_{R_1}$. Then on the annulus  $\{R_2'\leq |z|\leq R_1'\}$, there exists a constant $c>0$ such that 
\[-\frac{r+1-2l}{r}\log|q|\leq \widetilde v_l+c,\quad 1\leq l\leq n.\]

We construct the supersolution $\widetilde\vecw_{+}$ as follows:\\
on $|z|<R_2'$, let $\widetilde\vecw_{+}=\widetilde{\vecw}_{q}$;\\
on $ |z|\geq R_2'$, let $\widetilde\vecw_{+}=\min\{\widetilde{\vecw}_{q}, \widetilde\vecv+c\}.$
By the choice of $c$, $\widetilde\vecw_+=\widetilde\vecw_{q}$ in the neighborhood of $\partial B_{R_2'}$. So $\widetilde\vecw_{+}$ is continuous. Since $\widetilde\vecw_{q}, \widetilde\vecv+c$ are both supersolutions on $\{|z|\geq R_2'\}$, by Proposition \ref{Properties}, the minimum of $\widetilde\vecw_{q}$ and $\widetilde\vecv+c$ is a supersolution in the weak sense. Therefore $\widetilde\vecw_+\in  C^0(\mathbb C)\cap W^{1,2}_{loc}(\mathbb C)$ is a supersolution.

By Theorem \ref{CompletenessRealCurvature}, the solution $\widetilde\vecu$ satisfies $\widetilde\vecu<\widetilde\vecw_{q}$, since $|q|^{\frac{2}{r}}$ does not define a complete metric on $B_{R_1}$. It is then clear that $\widetilde\vecw_-<\widetilde \vecw_+$. Applying Proposition \ref{SuperSubExistence}, there is a $C^{\infty}$ solution $\widetilde\vecw$ satisfying $\widetilde\vecw_-\leq \widetilde\vecw\leq\widetilde \vecw_+$. Outside $B_{R_1}$, $\widetilde\vecv=\widetilde\vecw_-\leq \widetilde\vecw\leq \widetilde\vecw_+=\widetilde\vecv+c$. Since $\widetilde\vecv$ is a complete solution outside $B_{R_1}$,  the solution $\widetilde\vecw$ is complete on $\mathbb C$.
\hfill\qed 

\begin{prop}\label{ParabolicRiemannSurface}
Let $X$ be a parabolic Riemann surface with an Euclidean metric $g$ and a holomorphic $r$-differential $q$. If $q\neq 0$, then the system (\ref{systemhalf}) admits a unique complete solution. If $q=0$, then the system (\ref{eq;20.7.2.1}) admits no solution, that is, $\Toda(q,g)=\emptyset.$
\end{prop}
\pf
If $q=0$, the statement follows directly from Lemma \ref{NoSolutionOnParabolicSurface}. 

Suppose now $q\neq 0.$
Let $X$ be covered by $\mathbb C$ under the map $p:\mathbb C\rightarrow X$, with the covering transformation group of $X$ as $\Gamma<Aut(\mathbb C)=\{z\mapsto az+b, a\in \mathbb C^*,b\in \mathbb C\}$, i.e. $X=\mathbb C/ \Gamma$. Lift $q, g$ to $\hat q, \hat g$ on $\mathbb C$ which are invariant under $\Gamma.$ By Proposition \ref{ComplexPlaneExistence}, there exists a unique complete solution $\hat {\vecw}\in \Toda(\hat q, \hat g)$. By Proposition \ref{Pullback}, $\gamma^*(\hat\vecw)=\hat\vecw.$

Hence $\hat {\vecw}$ descends to a solution ${\vecw}$ on $\mathbb C/\Gamma=X$ which is still complete. 
\hfill\qed \\

Combining Theorem \ref{cor;20.9.18.110}, Proposition \ref{HyperbolicRiemannSurface} and Proposition \ref{ParabolicRiemannSurface}, we obtain
\begin{thm}\label{ExistenceUniquenessTheorem} Let $X$ be any non-compact Riemann surface with any K\"ahler metric
$g$ and a holomorphic $r$-differential $q$. We assume that $q\neq  0$ if $X$
is parabolic. Then, there exists a unique complete solution $\vecw\in
\Toda(q, g)$. Moreover, it is real.
 \hfill\qed 
 \end{thm}

In the case $q$ has finitely many zeros, we can construct  a solution whose asymptotic behavior is purely controlled by $q$ using the method of super-subsolution again. This solution is not necessarily complete. For example, when $X=\mathbb C$, $q=ze^z$, the solution is clearly not complete.
\begin{prop}\label{FiniteManyZerosExistence}
If $q$ has finitely many zeros, then there exists a smooth solution $\vecw$ of the system (\ref{systemhalf}) satisfying outside a relatively compact open set containing all zeros of $q$, 
$$\vecw_q-c\leq \vecw\leq \vecw_q.$$
\end{prop}
\pf
Denote by $Z$ the set of zeros of $q$. In the case $Z$ is empty, $\vecw_q$ is a solution as desired. 

In the case $Z$ is nonempty, we will use the method of super-subsolution for the system to show the existence by constructing a supersolution and a subsolution. 

For each point $P\in Z$, choose a neighbourhood $N_P$ of $P$ such that (i) $N_P$ is embedded into $\mathbb C$, which induces a coordinate $z_P$ on $N_P$, (ii) $\{|z_P|<1\}$ is relatively compact in $N_P$, (iii) $N_P\cap N_Q=\emptyset$ for $P, Q\in Z$. Let $\Omega_{P,2}=\{|z_P|<1\}.$ Choose neighbourhoods $\Omega_{P,0}$ and $\Omega_{P,1}$ of P such that $\Omega_{P, 0}\Subset\Omega_{P,1}\Subset\Omega_{P,2}$.
Let $\Omega_i=\cup_{P\in Z}\Omega_{P,i}(i=0,1,2)$. 
Since $Z$ is finite, $\Omega_i (i=0,1,2)$ are relatively compact. By Proposition \ref{HyperbolicRiemannSurface}, the system admits a complete solution $\vecu_{P}=(u_{P,1},\cdots, u_{P,n})$ on $\Omega_{P,2}$. Denote by $\vecu$ the solution on $\Omega_2$ whose restriction to $\Omega_{P,2}$ is $\vecu_P$.

The relatively compact subset $\Omega_1\setminus \Omega_0$ does not contain any zeros of $q$, then $-\frac{r+1-2l}{r}\log|q|_g (l=1,\cdots, n)$ is bounded from above. Therefore, there exists a constant $c>0$ such that on $\Omega_1\setminus\Omega_0$,
\begin{equation}\label{ChoiceOfc}
u_l\geq -\frac{r+1-2l}{r}\log|q|_g-c, \quad \forall 1\leq l\leq n.
\end{equation}

We construct a supersolution as follows: \\
on $\Omega_1$, let $\vecw_+=\min \{\vecu+c, \vecw_{q}\}$;\\
outside $\Omega_1$, let $\vecw_+=\vecw_{q}$. We must check the continuity of $\vecw_+$ on $\partial \Omega_1$ and at zeros of $q$. Because of the choice of $c$, $\vecu+c\geq \vecw_{q}$ in $\Omega_1\setminus \Omega_0$, so $\vecw_+=\vecw_{q}$  in the neighborhood of $\partial \Omega_1$. Since $\vecu$ is continuous on $\Omega_1$ and $\vecw_{q}$ goes to $+\infty$ at zeros of $q$, the vector function $\min\{\vecu+c,\vecw_{q}\}=\vecu+c$  is continuous at zeros of $q$. Since $\vecu$ is a solution in $\Omega_1$, $\vecu+c$ is a supersolution. By Proposition \ref{Properties}, on $\Omega_1$, $\vecw_+$ is still a supersolution. Hence on $X$, $\vecw_+$ is a supersolution.

We construct a subsolution as follows:\\
on $\Omega_0$, let $\vecw_-=\vecu$; \\
on the annulus $\Omega_2\setminus \Omega_0$, let $\vecw_-=\max \{\vecu, \vecw_{q}-c\}$; \\
outside $\Omega_2$, let $\vecw_-=\vecw_{q}-c$. We must check the continuity of $\vecw_-$ on $\partial \Omega_0$ and $\partial \Omega_2$. From the choice of $c$, $\vecu\geq \vecw_{q}-c$ in the neighborhood of $\partial \Omega_0$, so $\vecw_-=\vecu$ in the neighborhood of $\partial \Omega_0$. By Proposition \ref{HyperbolicRiemannSurface}, $-\frac{2}{r+1-2k}u_{P,k}(k=1,\cdots, n)$ and $\log (1-|z_P|^2)^{-2}$ are mutually bounded on $\Omega_{P,2}$. Hence $u_{P,k}(z_P)\rightarrow-\infty$ as $z_P\rightarrow\partial\Omega_{P,2}$. So $\vecw_-=\vecw_{q}-c$ in the neighborhood of $\partial\Omega_2$. Since $\vecw_{q}$ is a solution on $X\setminus Z$, $\vecw_{q}-c$ is a subsolution. By Proposition \ref{Properties}, on $\Omega_2$, $\vecw_-$ is still a subsolution. Hence on $X$, $\vecw_-$ is a subsolution.

By Theorem \ref{CompletenessRealCurvature}, $\vecu< \vecw_{q}$ holds on $\Omega_2$, since $Z$ is nonempty. It is then clear that $\vecw_-< \vecw_+$. Applying Proposition \ref{SuperSubExistence}, there is a $C^{\infty}$ solution $\vecw$ satisfying $\vecw_-\leq \vecw\leq \vecw_+$. The inequalities of $\vecw$ are clear.
\hfill\qed

\end{document}

%% file: notation.tex
\newcommand{\nbige}{\mathcal{E}}

\newcommand{\nbigo}{\mathcal{O}}
\newcommand{\nbigp}{\mathcal{P}}

\newcommand{\proj}{\mathbb{P}}
\newcommand{\seisuu}{{\mathbb Z}}

\newcommand{\cnum}{{\mathbb C}}
\newcommand{\real}{{\mathbb R}}

\newcommand{\hyperk}{\mathbb{K}}

\newcommand{\vecv}{{\boldsymbol v}}
\newcommand{\vecu}{{\boldsymbol u}}
\newcommand{\vecw}{{\boldsymbol w}}

\newcommand{\veca}{{\boldsymbol a}}

\newcommand{\lrarr}{\longrightarrow}

\newcommand{\pf}{{\bf Proof}\hspace{.1in}}
\newcommand{\qed}{\mbox{\rule{1.2mm}{3mm}}}

\def\End{\mathop{\rm End}\nolimits}

\def\Re{\mathop{\rm Re}\nolimits}

\def\Gr{\mathop{\rm Gr}\nolimits}

\def\rank{\mathop{\rm rank}\nolimits}

\def\Gr{\mathop{\rm Gr}\nolimits}

\def\tr{\mathop{\rm tr}\nolimits}
\def\Tr{\mathop{\rm Tr}\nolimits}

\def\id{\mathop{\rm id}\nolimits}

\newcommand{\del}{\partial}
\newcommand{\delbar}{\overline{\del}}

\newcommand{\baralpha}{\overline{\alpha}}
\newcommand{\alphabar}{\baralpha}

\def\Harm{\mathop{\rm Harm}\nolimits}

\def\rc{\mathop{\rm c}\nolimits}

\def\Toda{\mathop{\rm Toda}\nolimits}

\newcommand{\vbar}{\overline{v}}

\newcommand{\gtilde}{\widetilde{g}}

\newcommand{\wtilde}{\widetilde{w}}

\newcommand{\betabar}{\overline{\beta}}

\newcommand{\Xbar}{\overline{X}}



%% file: new_theorem.tex
\newtheorem{thm}{Theorem}[section]
\newtheorem{cor}[thm]{Corollary}

\newtheorem{rem}[thm]{Remark}
\newtheorem{lem}[thm]{Lemma}
\newtheorem{prop}[thm]{Proposition}
\newtheorem{df}[thm]{Definition}

\newtheorem{condition}[thm]{Condition}